\documentclass[12pt]{article}
\usepackage[english]{babel}
\usepackage{amssymb,amsmath,graphics}
\parskip\medskipamount
\newtheorem{theorem}{Theorem}[section]
\newtheorem{Theorem}{Main Result}
\newtheorem{remark}[theorem]{Remark}
\newtheorem{prop}[theorem]{Proposition}
\newtheorem{lemma}[theorem]{Lemma}
\newtheorem{cor}[theorem]{Corollary}
\def\<{\langle}
\def\>{\rangle}
\newcommand{\proof}{\emph{Proof.~}}
\newcommand{\dd}{\mathsf{d}}

\newcommand{\cF}{\mathcal{F}}

\def\qed{{\hfill\hphantom{.}\nobreak\hfill$\Box$}}

\newcommand{\A}{\mathbf{A}}
\newcommand{\R}{\mathbb{R}}

\newcommand{\inc}{\mbox{\tt I}}

\newcommand{\opp}{\mathrm{opp}}
\begin{document}

\author{Koen Struyve\thanks{The first author is supported by  the Fund for Scientific Research --
Flanders (FWO - Vlaanderen)} \and Hendrik Van Maldeghem\thanks{The second author is partly supported by a
Research Grant of the Fund for Scientific Research -- Flanders (FWO - Vlaanderen)} }
\title{\bf Two-dimensional affine $\R$-buildings defined by generalized polygons with non-discrete valuation}

\maketitle

\begin{abstract} In this paper we complete the proof of the `equivalence' of non-discrete $\R$-buildings 
of types $\widetilde{\mathsf{A}}_2$ and $\widetilde{\mathsf{C}}_2$, with, respectively,  projective planes and generalized quadrangles with non-discrete valuation, begun in \cite{kshvm}. We also complete the proof of the `equivalence' of an affine building of rank 3 with a generalized polygon with discrete valuation (by proving this for generalized hexagons), begun in \cite{Mal:90}. We also complement the main result of \cite{Mal:92} by proving uniqueness up to scalar multiples of the weight sequences of polygons with non-discrete valuation. As an application, we produce some new explicitly defined non-discrete $\R$-buildings, in particular a class of type $\widetilde{\mathsf{A}}_2$ with arbitrary residues.
\end{abstract}

\section{Introduction}

In 1984, Jacques Tits \cite{tits84} classified affine buildings of rank at least 4. In fact, he also included in his work the so-called \emph{non-discrete affine buildings}, which he called \emph{syst\`emes d'appartements}, or \emph{apartment systems}. Basically, these are building-like structures with one big difference: they are no longer simplicial. Easy examples are $\R$-trees (rank 2 case; these are trees that continuously branch), or the ``buildings'' related to the ``parahoric'' subgroups of a Chevalley group over a field with non-discrete valuation. From the geometric point of view, the case of rank 3 --- when the apartments are 2-dimensional --- is very interesting since non-classical phenomena occur there. 

In \cite{tits84} Tits associates to every \emph{symmetric} apartment system a so-called \emph{building at infinity}, which is a simplicial spherical building, see also \cite{Bru-Tit:72}. The rank of this building at infinity is precisely the dimension of its apartments. Hence, in the 2-dimensional case, generalized polygons appear. When the apartment system is irreducible, then this polygon is not a digon. In the simplicial case, the only generalized polygons that occur are projective planes, generalized quadrangles and generalized hexagons. In 1992, the second author \cite{Mal:92} introduced the notion of a generalized polygon with discrete valuation and conjectured that the resulting polygons are precisely the buildings at infinity of the rank 3 irreducible affine buildings. This conjecture was verified for all cases except for the generalized hexagons.   However, in \cite{kshvm}, we already showed that any generalized hexagon isomorphic to the building at infinity of an affine building, admits a discrete valuation in the sense of \cite{Mal:92}. The starting point of the present paper is to complete the proof of this conjecture by showing that every generalized hexagon with valuation is isomorphic to the building at infinity of an affine building (of type $\widetilde{\mathsf{G}}_2$).  

But we achieve more. In \cite{Ber-Kap:08}, Berenstein and Kapovich prove the existence of 2-dimensional (nontrivial, i.e., no blow-ups of spherical buildings) apartment systems 
admitting a generalized $n$-gon at infinity for any integer value of $n>2$. The natural question hereby is whether these structures are also characterized by admitting a valuation in some sense. Notice that discrete valuations are nonexistent for $n$-gons with $n\neq 3,4,6$ by \cite{Mal:92}.  However, as we shall show below, if we symmetrize the definition of generalized polygon with valuation (with respect to the notions of points and lines) and allow real values (we shall call these \emph{generalized $n$-gons with real valuation}), then the only weight sequences (for a definition see below) that can occur are the ones that come from 2-dimensional apartment systems as shown in \cite{kshvm}. Moreover, if $n=3,4$, then we provide a detailed proof for the complete equivalence between generalized $n$-gons with real valuation and 2-dimensional symmetric affine apartment systems. As an application we construct classes of explicit examples of such structures which are not of Bruhat-Tits type, and which include locally finite ones. These constructions are similar to the constructions due to the second author in the simplicial case, see \cite{Mal:87,Mal:88,Mal:89,Mal:90}. 

Remarkably, as a byproduct, we obtain that projective planes with real valuation are equivalent with ultra-metric planes in which  all triangles satisfy the sine rule, for an appropriate though natural definition for angles between lines. 

In the ideal case, one would like to prove the conjecture that the just mentioned equivalence holds for all $n\geq 3$. However, this seems to be out of reach for now. In our present approach, the complications in the proofs seem to grow exponentially with the girth. For $n=5$, it is just feasible, but too long to include here. For $n=6$, assuming discreteness allows for an alternative argument, as we shall see. Notice that our proofs for $n=3,4$ provide different arguments for the simplicial case, which are in fact drastically shorter and more direct than the original proofs of the second author. One does not need to go around the \emph{Hjelmslev geometries} and the rather complicated axiomatization related to this (see e.g.~\cite{Han-Mal:90}). These geometries were needed to define the vertices of the affine building. In the present approach, we do no longer have vertices, but the points of the apartment system are the different valuations that emerge from the given one. This simple idea, however, requires a lot of unavoidable technicalities to take care of.  For example, it is already fairly technical to prove that the residue of an $n$-gon with valuation is again a generalized $n$-gon. We will do this explicitly for $n\leq 6$. It will be clear that similar methods should work in general, but our present approach fails for that. So, on the one hand, the present methods are significantly stronger than the old ones developed by the second author in the eighties, on the other hand, one needs an improvement of another magnitude to prove the full conjecture. 

Finally, we would like to remark that the central objects in this paper are inventions of Jacques Tits, without whom this paper would never have been written. The classification of irreducible affine buildings of rank at least 4 was just completed when the second author started a PhD, greatly inspired by this, on affine buildings of type    $\widetilde{\mathsf{A}}_2$, advised by Mark Ronan in Chicago. The second author also wants to express his profound thanks and estimation for the work and especially the lectures of Jacques Tits at the Coll\`ege de France. The latter were a constant motivation and inspiration. How a whole career can be `build' on investigating left-overs of one of the most remarkable mathematicians of the last century, the one that created Incidence Geometry and made it almost a synonym of Group Theory. 

\section{Preliminaries and Main Results}\label{MR}
\subsection{$\R$-buildings}
Let $(\overline{W},S)$ be a finite irreducible Coxeter system. So $\overline{W}$ is presented by the set $S$ of
involutions subject to the relations which specify the order of the products of every pair of involutions. This
group has a natural action on a real vector space $V$ of dimension $|S|$. Let $\A$ be the affine space associated to
$V$. We define $W$ to be the group generated by the translations of $\mathbf{A}$ and $\overline{W}$.

Let $\mathcal{H}_0$ be the set of hyperplanes of $V$ corresponding to the axes of the reflections in $S$ and
all its conjugates. Let $\mathcal{H}$ be the set of all translates of all elements of $\mathcal{H}_0$. The
elements of $\mathcal{H}$ are called \emph{walls} and the (closed) half spaces they bound are called
\emph{half-apartments} or \emph{roots}. A \emph{vector sector} is the intersection of all roots that (1) are
bounded by elements of $\mathcal{H}_0$, and (2) contain a given point $x$ that does not belong to any element
of $\mathcal{H}_0$. The bounding walls of these roots will be referred to as the \emph{side-walls} of the
vector sector. A vector sector can also be defined as the topological closure of a connected component of
$V\setminus(\cup\mathcal{H}_0)$. Any translate of a vector sector is a \emph{sector}, with corresponding
translated \emph{side-walls}. A \emph{sector-facet} is an infinite intersection of a given sector with a finite
number of its side-walls. This number can be zero, in which case the sector-facet is the sector itself; if this
number is one, then we call the sector-facet a \emph{sector-panel}. The intersection of a sector with all its
side-walls is a point which is called the \emph{source} of the sector, and of every sector-facet defined from
it.

An \emph{$\R$-building} (also called an \emph{affine apartment system}) (definition by Jacques Tits as can be
found in~\cite{ronan} by Mark Ronan, along with some historic background) is an object $(\Lambda,\cF)$
consisting of a set $\Lambda$ together with a collection $\cF$ of injections of $\mathbf{A}$ into $\Lambda$ (called \emph{charts})
obeying the five conditions below. The image of $\mathbf{A}$ under an $f\in \cF$ will be called an
\emph{apartment}, and the image of a sector, half-apartment, $\dots$ of $\mathbf{A}$ under a certain $f\in \cF$
will be called a \emph{sector}, \emph{half-apartment}, $\dots$ of $\Lambda$.
\begin{itemize}
 \item[(A1)] If $w\in W$ and $f\in \cF$, then $f \circ w \in \cF$.
 \item[(A2)] If $f,f' \in \cF$, then $X=f^{-1}(f'(\mathbf{A}))$ is closed and convex
 in $\textbf{A}$, and $f|_X = f'\circ w|_X$ for some $w\in W$.
 \item[(A3)] Any two points of $\Lambda$ lie in a common apartment.
\end{itemize}
The last two axioms allow us to define a function $\dd : \Lambda \times \Lambda \rightarrow \R^+$ such that for
any $a,b \in \A$ and $f\in \cF$, $\dd(f(a),f(b))$ is equal to the Euclidean distance between $a$ and $b$ in
$\A$.
\begin{itemize}
 \item[(A4)] Any two sectors contain subsectors lying in a common apartment.
 \item[(A5$'$)] Given $f \in \cF$ and a point $\alpha \in \Lambda$, there is a retraction
 $\rho : \Lambda \rightarrow f(\A)$ such that the preimage of $\alpha$ is $\{\alpha\}$ and which diminishes $\dd$.
\end{itemize}

We call $|S|$, which is also equal to $\dim \mathbf{A}$, the \emph{dimension} of $(\Lambda,\cF)$. We will
usually denote $(\Lambda,\cF)$ briefly by $\Lambda$, with slight abuse of notation.

A detailed analysis of this definition and variations of it has been carried out by Anne Parreau in~\cite{parreau}.
In particular, she shows that, if Conditions (A1), (A2), (A3) and (A4) are satisfied, then (A5$'$) is equivalent to  $\dd$ being a distance function, together with
\begin{itemize}
\item[(A5)] If we have three apartements, such that each two apartments of these share a half-apartment, then the intersection of all three is non-empty. 
\end{itemize}

\subsection{Generalized polygons}

Generalized polygons are the geometries corresponding to the spherical rank 2 buildings. Since we will use some
specific terminology of these geometries, we introduce this now.

A \emph{generalized $n$-gon}, $n\in\mathbb{N}$, $n\ge 2$, or \emph{generalized polygon} $\Gamma=(P,L,\inc)$ is
a structure consisting of a \emph{point set} $P$, a \emph{line set} $L$ (with $P\cap L=\emptyset$), and a
symmetric \emph{incidence relation} $\inc$ between $P$ and $L$, turning $P\cup L$ into a bipartite graph
$\mathfrak{G}$ satisfying the following axioms.
\begin{itemize}
 \item[(GP1)] Every element is incident with at least three other elements.
 \item[(GP2)] For every pair of elements $x,y\in P\cup L$,
 there exists a sequence $x_0=x,x_1,\ldots,x_{k-1},x_k=y$, with $x_{i-1}\inc x_{i}$ for $1\leq i\leq k$ and
 with $k\leq n$.
 \item[(GP3)] The sequence in (GP2) is unique whenever $k<n$.
\end{itemize}

If instead a weaker version of (GP1) is satisfied where each element is incident with at least two elements, we speak about a \emph{weak generalized $n$-gon}.

A \emph{path} is an ordered set of elements such that each two subsequent elements in the set are incident. The \emph{length of a path} is the number of elements in the set minus one. A path is \emph{closed} if the last element of the set equals the first, and is \emph{non-stammering} if for each element of the ordered set, the two neighbours are different.

The \emph{distance} $\dd(x,y)$ between two elements $x,y$ is the length of a shortest path between both. If two elements $x$ and $y$ are at distance 2, then $xy$ will denote the unique element incident with both.
 
Two points are \emph{collinear} if they are incident with a common line, two lines are \emph{concurrent} if they are incident with a common point. Two elements are \emph{adjacent} if they are collinear or concurrent.

If two elements are at distance $n$, they are called \emph{opposite}. If two elements $x$ and $y$ are not opposite, then the unique element incident with $y$ closest to $x$ is the \emph{projection of $x$ on $y$}.

\subsection{Generalized polygons with (non-discrete) valuation}
Let $\Gamma=(P,B,\inc)$ be a generaized $n$-gon and $u$ a function called the \emph{valuation} acting on pairs of adjacent elements, and with images in $\R^+ \cup \{ \infty \} $ ($\R^+$ being the non-negative real numbers, and using the natural order on this set with $\infty$ as largest element). Then we call $(\Gamma,u)$ an \emph{$n$-gon with (non-discrete) valuation} and \emph{weight sequence} $(a_1,a_2, \dots, a_{n-1}, a_{n+1},a_{n+2}, \dots a_{2n-1} ) \in (\R^+_0)^{2n-2}$ ($\R^+$ being the positive real numbers) if the following conditions are met:
\begin{itemize}
\item[(U1)] For each element $z$, there exists a pair $x$ and $y$ of elements incident with $z$ such that $u(x,y)=0$.
\item[(U2)] $u(x,y) = \infty$ if and only if $x=y$.
\item[(U3)] If $x,y$ and $z$ are collinear points or concurrent lines, then $u(x,y) < u(y,z)$ implies $u(x,z) = u(x,y)$.
\item[(U4)] Whenever $x_1 \inc x_2 \inc \dots \inc  x_{2n} \inc x_1$, with $x_i \in P \cup B$, one has 
\begin{equation}
\sum_{i=1}^{n-1} a_i u(x_{i-1},x_{i+1}) = \sum_{i=n+1}^{2n-1} a_i u(x_{i-1},x_{i+1})
\end{equation}
\end{itemize}
One direct implication of (U2) and (U3) is that $u$ is symmetric (by putting $x=z$ in (U3)). Also remark that this definition is self-dual when interchanging lines and points, so whenever a statement is proven, we also have proven the dual statement.

If we speak about the valuation of a side or corner $x$ in an ordinary $n$-gon $\Omega$ we mean the valuation between respectively the two corners or sides incident with $x$ in $\Omega$. If we talk about the valuations in an ordinary $n$-gon, then we mean all the valuations of sides and corners. A path $(x_0,x_1, \dots, x_m)$ is said to have valuation zero if $u(x_{i-1},x_{i+1}) = 0$ for each $i \in \{1,2\dots, m-1\}$. Such a path is also non-stammering. We now show some preliminary lemmas that we will use to formulate one of the main results.
\begin{lemma}\label{lemma:pq0}
Given a line $L$ and a point $p \inc L$, then there exists a point $q \inc L$ such that $u(p,q) = 0$.
\end{lemma}
\proof
Due to (U1) there exist two points $r,s \inc L$ such that $u(r,s) = 0$. Applying (U3) we obtain that either $u(p,r)= 0$ or $u(p,s) =0$, in each case we have found a suitable $q$.
\qed

\begin{lemma}\label{lemma:contain}
Each path $(x_0,x_1, \dots ,x_m)$ with $m \leq n+1$ and valuation zero is contained in an ordinary $n$-gon $\Omega$ where all the valuations of corners and sides are zero.
\end{lemma}
\proof
Using the previous lemma we can extend the path to a path $(x_0 := p, x_1 := L, \dots, x_n,x_{n+1})$ with valuation zero. It is now easily seen that the other valuations in the unique ordinary $n$-gon triangle spanned by the path are zero too by (U4).
\qed

In order to make notations easier, such ordinary $n$-gon with all valuations zero will be referred to as a \emph{non-folded $n$-gon}. If there are exactly two non-zero valuations in (necessarily) opposite elements $x$ and $y$ of an ordinary $n$-gon, then this ordinary $n$-gon will be referred to as a \emph{simply folded $n$-gon folded along $x$ (or $y$)}, two elements in such an $n$-gon at the same distance from $x$ (and hence also at the same distance from $y$) are said to be \emph{folded together} in that $n$-gon. The first main result will imply that $a_1=a_{n+1}$ and thus that the valuations in $x$ and $y$ are equal due to (U4). 

Two opposite elements in $\Gamma$ are said to be \emph{residually opposite} if there is a shortest path between them with valuation zero. If this is the case, then by (U4) all shortest paths between both elements have valuation zero. If $x$ is an an element of $\Gamma$ then $[x]_\opp$ is the set of residually opposite elements, this set is non-empty due to the previous lemma. We say that two elements $x$ and $y$ as \emph{residually equivalent} if $[x]_\opp = [y]_\opp$. The equivalence class is denoted as $[x]=[y]$. It is clear that all elements of one equivalence class share the same type, so these classes can be referred to as \emph{residual points} ($[P]$) or \emph{residual lines} ($[B]$) depending on this type. A residual point $[p]$ is said to be incident with a residual line $[L]$ if there are $p' \in [p]$ and $L' \in [L]$ such that $p' \inc L'$. We then write $[p] \inc_r [L]$. The geometry $\Gamma_r([P],[B],I_r)$ is denoted as the \emph{residue} defined by $u$. The distance $\dd_r$ in the incidence graph of this geometry is called the \emph{residual distance}.

\subsection{Main Results}
Let $(\Gamma,u)$ be a generalized $n$-gon with (non-discrete) valuation and weight sequence $(a_1,a_2, \dots, a_{n-1}, a_{n+1},a_{n+2}, \dots a_{2n-1} )$. 
\begin{Theorem}
If $u$ has non-zero values, the weight sequence $(a_1,a_2, \dots, a_{n-1}, a_{n+1},a_{n+2}, \dots a_{2n-1} )$ is a multiple of the weight sequence $(b_1,b_2, \dots, b_{n-1}, b_{n+1},b_{n+2}, \dots b_{2n-1} )$ with $b_i = |\sin(i \pi / n) |$. 
\end{Theorem}

\begin{Theorem}
If $3 \leq n \leq 6$, then the residue defined by $u$  is a (weak) generalized $n$-gon.
\end{Theorem}

\begin{Theorem}
If $n\in\{3,4\}$, or if $n=6$ and $u$ is discrete, then there exists a two-dimensional $\R$-building $(\Lambda,\cF)$ such that $\Gamma$ is isomorphic to the generalized polygon at infinity of $(\Lambda,\cF)$ with valuation as defined in~{\rm \cite{kshvm}}.
\end{Theorem}

\subsection{An application to ultrametric projective planes}\label{section:app}
In this application we explore a surprising link between projective planes with valuations and some geometric conditions from Euclidean geometry.

Suppose $(\Gamma,u)$ is a trigon (or projective plane) with valuation. Choose $t \in \R$ with $t>1$. We then can define a function $\dd(p,q) = t^{-u(p,q)} \in [0,1]$ on pairs of points, and a similar function $\angle(L,M) = \arcsin(t^{-u(L,M)}) \in [0,\pi /2]$ on pairs of lines.

\begin{theorem}\label{theorem:app}
A projective plane $\Gamma$ with a distance function $d$ on pairs of points valued in $[0,1]$ and an angle function $\angle$ on pairs of lines valued in $[0,\pi /2]$, is constructed from a projective plane with valuation as above, and hence is isomorphic to the building at infinity of some $\R$-building, if and only if the following conditions are fulfilled.
\begin{itemize}
\item[{\rm (M1)}] $\dd$ is an ultrametric (this is a metric satisfying the stronger triangular inequality $\dd(p,q) \leq \max(\dd(p,r), \dd(r,q))$).
\item[{\rm (M2)}] Two lines have angle zero if and only if they are equal.
\item[{\rm (M3)}] On each line there are two points on the maximal distance 1 from each other.
\item[{\rm (M4)}] Through each point there are two lines with a right ($\pi /2$) angle.
\item[{\rm (M5)}] The sine rule is fulfilled, i.e. if we have a triangle with lengths of the sides $A$, $B$ and $C$ and opposing angles $\alpha$, $\beta$ and $\gamma$, then 
\begin{equation*}
\frac{A}{\sin \alpha}=\frac{B}{\sin \beta}=\frac{C}{\sin \gamma}.
\end{equation*}
\end{itemize}
\end{theorem}

\section{Proof of the first main result}
We start with a polygon $\Gamma$ with valuation $u$, with weight sequence $(a_1,a_2, \dots, a_{n-1}, a_{n+1},a_{n+2}, \dots a_{2n-1} )$, and such that $u$ has non-zero values. Our proof is heavily inspired by a similar result for the discrete case in~\cite{Mal:92}. In fact, we will use some of the results (with the proofs remaining valid in the non-discrete case) obtained there, directly in our proof. In particular, and to begin with, it is shown in 3.1 of \cite{Mal:92} that the weight sequence is unique, up to a non-zero multiple. As is also exploited in \cite{Mal:92}, this has as consequence that the weight sequence is symmetric, i.e., $a_i= a_{n-i} =a_{n+i}=a_{2n-i}$ for $i \in \{1,2, \dots, n-1\}$.

Now let $(x_0,x_1, \dots, x_{2n} = x_0)$ be any closed path of length $2n$ in $\Gamma$. Because of (U4) we know that  
\begin{equation}
\sum_{i=1}^{n-1} a_i u(x_{i-1},x_{i+1}) = \sum_{i=n+1}^{2n-1} a_i u(x_{i-1},x_{i+1}),
\end{equation}
and also that
\begin{equation}
\sum_{i=3}^{n+1} a_{i-2} u(x_{i-1},x_{i+1}) = \sum_{i=n+3}^{2n+1} a_{i-2} u(x_{i-1},x_{i+1}).
\end{equation}
If one takes the sum of both equations, and simplifies the resulting expression using $a_1= a_{n-1}=a_{n+1}=a_{2n-1}$, one obtains
\begin{align}
\begin{split}
a_2 u(x_1,x_3) + \sum_{i=3}^{n-1} (a_i + a_{i-2}) u(x_{i-1},x_{i+1}) + a_{n-2}u(x_{n-1},x_{n+1}) \\ =  a_{n+2} u(x_{n+1},x_{n+3}) +  \sum_{i=n+3}^{2n-1} (a_i + a_{i-2}) u(x_{i-1},x_{i+1}) + a_{2n-2} u(x_{2n-1},x_{2n+1}).
\end{split}
\end{align}
This implies that $$(a_2,a_3+a_1,a_4+a_2, \dots, a_{n-1}+a_{n-3},a_{n-2}, a_{n+2},a_{n+3}+a_{n+1}, \dots, a_{2n-1}+a_{2n-3},a_{2n-2})$$ is also a weight sequence. Hence there exists some positive real number $k$ satisfying

\begin{equation}\label{eq5}
\left\{ \begin{array}{l}
ka_1 = a_2, \\
ka_2 = a_3+a_1,\\
ka_3 =a_4+a_2, \\
\dots \\
ka_{n-2} =a_{n-1}+a_{n-3}, \\
ka_{n-1} = a_{n-2}.\\
\end{array} \right.
\end{equation}
One notices that, by taking the sum of all equations in the system of equations above, that
\begin{equation}
k \sum_{i=1}^{n-1} a_{i} = 2 \sum_{i=1}^{n-1} a_{i} - (a_1 +a_{n-1}).
\end{equation}
This implies that $1\leq k < 2$. As a consequence, we can find an $\alpha \in ]0, \pi/3]$ such that $k=2 \cos \alpha$. Also remark that $a_j =k a_{j-1} - a_{j-2}$ for $j \in \{3,n-1\}$. If we formally set $a_0=a_n=0$, then this is also true for $j \in \{2,n\}$. Furthermore we can suppose that $a_1 = \sin \alpha$.

\begin{lemma}
For $i \in \{0,1,\dots, n\}$ we have $a_i = \sin(i \alpha)$. 
\end{lemma}
\proof
We prove this using induction on $i$. It is clear that this holds for $i=0$ and $i=1$ (by assumption and by definition of $\alpha$, respectively). So let $i \geq 2$ such that $a_j = \sin ja$ for $j <i$. Then we know that:
\begin{align}
a_i =& ka_{i-1} - a_{i-2} \\
=& 2\cos \alpha \sin[(i-1)\alpha] - \sin[(i-2)\alpha] \\
=& \sin i\alpha
 \end{align}
The second equality follows from the induction hypothesis, the third from the trigonometric formula $\sin a + \sin b= 2 \sin [(a+b)/2] \cos [(a-b)/2]$. \qed

\begin{lemma}
$\alpha = \pi/n$.
\end{lemma}
\proof
We have that $a_n =0$, so $\sin n \alpha =0$ by the previous lemma. This yields $\alpha =m \pi /n$, with $m \in \mathbb{N}_0$ smaller than or equal to $n/3$ (since $\alpha \in ]0,\pi/3]$). At the same time we have $a_i >0$ for $i \in \{1, \dots, n-1 \}$. Let $t$ be the smallest integer greater than or equal to $n/m$. Because $n/m \leq t \leq 2n/m$ (by $n/m \geq 3$), it holds that $t m \pi/n \in [\pi, 2\pi]$, so $a_t \geq 0$. As $t$ clearly is in $\{1,2,\dots,n\}$, we obtain that $t=n$, which implies that $m=1$ (because $m\in \mathbb{N}_0$ and $n\geq 3$) and $\alpha = \pi/n$.
\qed  

Combining the two previous lemmas, we obtain: 

\begin{cor}
For $i \in \{0,1,\dots, n\}$: $a_i = \sin(i \pi /n)$, and any other weight sequence of $(\Gamma,u)$ is a mulptiple hereof.
\end{cor}

\begin{remark} \rm It is easy to see that all $k\in\R$ satisfying Equation~\ref{eq5} are precisely the eigenvalues of the path graph $P_{n-1}$ of length $n-2$, consisting of $n-1$ vertices.  Moreover, since all $a_i$ are positive, it is the unique eigenvalue for which the coordinates of the associated eigenvectors have constant sign. This observation can be used to give an alternative proof of the previous corollary. Doing so, one sees that $2\cos(\pi/n)$ is in fact the largest eigenvalue of $P_{n-1}$. 
\end{remark}

\section{Proof of the second main result}
By the first main result we can suppose for the proof of the second and third main result that the weight sequence is given by $a_i = |\sin(i \pi /n)| / \sin(\pi / n)$. In particular, $a_1=1$.

Let $n$ be a natural number with $3\leq n\leq 6$ for the rest of the proof.

If $x$ and $y$ are opposite elements, let $\tau(x,y)$ be the sum $\sum_{i=1}^{n-1} a_i u(x_{i-1},x_{i+1})$ where $(x_0 = x,x_1, \dots, x_{n-1},x_n = y)$ is a shortest path from $x$ to $y$; (U4) guarantees independence of the chosen path. 


Two elements $x$ and $y$ are said to be \emph{$t$-residually equivalent}, if for each element $z$ the following are equivalent: \begin{itemize} \item $z$ is opposite $x$ and $\tau(x,z) < t$; \item $z$ is opposite $y$ and $\tau(y,z)<t$. \end{itemize}

Notice that when $t=0$, this definition is trivially fulfilled.  

\begin{lemma}\label{lemma:zero}
Two adjacent elements $x$ and $y$ are $u(x,y)$-residually equivalent, but not $t$-residually equivalent with $t >  u(x,y)$.
\end{lemma}
\proof
Let $z$ be an element opposite $x$ with $\tau(x,z) < u(x,y)$. Consider the unique shortest path $(x_0=x,x_1=xy,x_2,\dots,x_n=z)$ from $x$ to $z$ containing $xy$. Because $a_1=1$, it holds that $u(x,x_2) \leq \tau(x,z) < u(x,y)$, so $u(y,x_2) = u(x,x_2)$ by (U3). This implies that $y$ and $z$ are opposite and that $\tau(y,z)=\tau(x,z)$ (the last is easily seen when considering the path $(y,x_1,x_2,\dots,x_n=z)$). 

If $t > u(x,y)$, then consider a path $(x,xy, y=y_2, \dots, y_n)$ where the path $(y_2, \dots, y_n)$ has valuation zero (possible by Lemma~\ref{lemma:pq0}). 
\qed

\begin{cor}\label{cor:zero}
If $x\inc y \inc z$, then $[x]=[z]$ if and only if $u(x,z) > 0$ .
\end{cor}

\begin{lemma}\label{lemma:path}
Given a closed path $\Psi$, then there are at least two sides having the same minimal valuation among all sides in $\Psi$.
\end{lemma}
\proof 
Let $x$ and $y$ be the two points on a side with minimal valuation, and suppose all other sides have valuation strictly larger than $u(x,y)$. Let $t$ be the second smallest valuation among the sides in $\Psi$. By repeatedly using Lemma~\ref{lemma:zero} and going from $x$ to $y$ in $\Psi$ not using $xy$, one proves that $x$ and $y$ are $t$-residually equivalent, which contradicts Lemma~\ref{lemma:zero}.
\qed


\begin{lemma}
If two elements $x$ and $y$ are not residually equivalent, but if there exist $a \inc x$ and $b \inc y$ which are residually equivalent, then there is an element $z$ residually opposite one element of $\{x,y\}$, but at distance $n-2$ from the other. 
\end{lemma}
\proof
Without loss of generality, one can suppose that there exists an element $d$ which is residually opposite $x$, but not residually opposite $y$. 

According to Lemma~\ref{lemma:pq0}, there exists an element $c$ incident with $x$ such that $u(a,c)=0$.  Let $(x=x_0, c= x_1, \dots, x_{n-1},d= x_n)$ be the unique shortest path from $x$ to $d$ containing $c$. The element $x_{n-1}$ is residually opposite, and thus also opposite, $a$ and $b$. This implies that $\dd(y,d)=n$ or $\dd(y,d)=n-2$. In the second case we are done, so suppose we are in the first case. Let $(y=y_0, y_1, \dots, y_{n-2}, y_{n-1}= x_{n-1}, y_n = d)$ be the unique shortest path from $y$ to $d$ containing $x_{n-1}$. Because the element $x_{n-1}$ is residually opposite $b$, the path $(b, y=y_0, y_1, \dots, y_{n-2}, y_{n-1}= x_{n-1})$ has valuation zero.  As $y$ is not residually opposite $d$, the valuation $u(y_{n-2},d)$ has to be non zero. So $x_{n-2} \neq y_{n-2}$ and $u(x_{n-2},y_{n-2}) =0$. The element $x_{n-2}$ will now be the desired element $z$, because it is residually opposite $y$, but at distance $n-2$ from $x$. \qed

\begin{lemma}
Let $\Omega$ be a simply folded $n$-gon. If two elements $x$ and $y$ are folded together in $\Omega$, then they are residually equivalent.
\end{lemma}
\proof
Here we need to distinguish between the different possibilities for $n$. Let $z$ be an element of $\Omega$ such that $\Omega$ is folded along $z$. 
\begin{itemize}
\item
$n=3$. For this case the result follows directly from Corollary~\ref{cor:zero}.
\item
$n=4$. Again using Corollary~\ref{cor:zero}, one only needs to prove that the two elements of $\Omega$ at distance $2$ from $z$ are residually equivalent. Suppose this is not the case. Using the previous lemma, one can assume without loss of generality that there is an element $a$ residually opposite $x$, but at distance 2 from $y$. 

Let $(x, xz,x_2,x_3,a)$ be the unique shortest path (which has valuation zero) from $x$ to $a$ containing $xz$. Let $z'$ be the element opposite $z$ in $\Omega$. The element $x_3$ is residually opposite $xz'$, and thus also residually opposite $yz'$ due to Corollary~\ref{cor:zero}. This implies that the valuations $u(y,a)$ and $u(x_3,ay)$ are zero. But as also the valuations $u(xz,x_3)$ and $u(x_2,a)$ are zero, (U4) would imply that $u(xz,zy)=0$, which is a contradiction.
\item
$n=5$. Using Corollary~\ref{cor:zero} and the previous lemma, one can assume without loss of generality that $x$ and $y$ are at distance $2$ from $z$, and that there exists an element $a$ residually opposite $x$, but at distance 3 from $y$.  

Let $(x, xz,x_2,x_3,x_4, a)$ be the unique shortest path (which has valuation zero) from $x$ to $a$ containing $xz$, and let $(y, y_1, y_2,a)$ be the shortest path from $y$ to $a$. Choose an element $b \inc a$ such that $u(b,x_4) =0$ (this is possible due to Lemma~\ref{lemma:pq0}). The element $xz$ is residually opposite $b$, thus so is $yz$. All of this implies that the path $(yz,y,y_1,y_2,a,b)$ has valuation zero.  A consequence is that $u(x_4,y_2) >0$, otherwise we could have chosen $b$ to be $y_2$, leading to a contradiction.

Let $z'$ be the element opposite $z$ in $\Omega$, and let $x', y'$ be the elements incident with $z'$ closest to $x$ and $y$ respectively. Now $x'$ and $y'$ are both residually opposite $x_3$, implying that the unique shortest path from $yy'$ to $x_3$ has valuation zero. If we look in the unique ordinary pentagon containing $yy', x_3$ and $y_2$, we see that the valuation of $x_3$ in this pentagon is non-zero because of (U4) and $u(x_4,y_2)>0$. By (U3) we then obtain that the valuation of $x_3$ in the unique ordinary pentagon containing $x_3$, $yy'$ and $z$ is zero. This contradicts (U4) and the fact that the valuation of $z$ in this pentagon is non-zero.

\item
$n=6$. Apart from the case handled in Corollary~\ref{cor:zero}, there are two cases to consider here.
\begin{itemize}
\item
The first case is when $x$ and $y$ lie at distance 2 or 4 from $z$, without loss of generality one can suppose this to be 2. Similarly to the previous cases, let $a$ be an element residually opposite $x$, but at distance 4 from $y$. Let $x_1$ be the unique element of $\Omega$ at distance 1 from $x$ and 3 from $z$. Now consider the unique shortest path $(x,x_1,x_2,x_3,x_4,x_5,a)$ from $x$ to $a$ containing $x_1$, and the unique shortest path $(y,y_1,y_2,y_3,a)$ from $y$ to $a$. Observe that $x_4 \in [z]_\opp$. Let $\Omega'$ be the unique ordinary simply folded hexagon containing $z$, $x_4$, $x$ and $yz$, and let $b$ be the element opposite $x_2$ in this hexagon. By (U3), the unique ordinary hexagon containing $y$, $b$, $y_1$, and $x_4$ is non-folded, so $u(y,b)$ is zero and $x_4 \in [y]_\opp$. 

Let $\Omega''$ be the unique ordinary hexagon containing $z$, $y$ and $x_3$, and $\Omega'''$ the unique ordinary hexagon containing $y$, $b$ and $x_3$. Let $c$ and $c'$ respectively be the elements opposite $xz$ in the hexagons $\Omega$ and $\Omega'' $ respectively. Let $d$ and $d'$ be the projections of $c$ and $c'$, respectively, on $y$. The hexagon $\Omega'''$ is a simply folded hexagon folded along $y$ (remember that $u(y,b)$ was zero). So $u(yz,d')$ is non-zero, and thus $u(d,d')$ is zero. This implies that $c \in [c']_\opp$, so also the element $c''$ opposite $yz$ in $\Omega$ is in $[c']_\opp$. Because the unique path from $c''$ to $c'$ containing $x_2$ has valuation zero, also the path from $xz$ to $c'$ containing $x$ has valuation zero. Thus $xz \in [c']_\opp$, which gives $yz \in [c']_\opp$ which is a contradiction because $yz$ and $c'$ are at distance 4 from each other.
\item
The last case to handle is the case where $x$ and $y$ are at distance 3 from $z$. For the final time, consider an element $a\in [x]_\opp$ and at distance 4 from $y$. Let $x'$ and $y'$ be the projections from $z$ on $x$ and $y$, respectively, and let $x''$ and $y''$ be the elements in $\Omega$ at distance 4 from $z$ and 1 from $x$ and $y$, respectively. Let $a'$ be the projection of $x''$ on $a$; this element is residually opposite $x'$, so it is also residually opposite $y'$ (as shown in the previous case). The unique shortest path from $y'$ to $a'$ containing $a$ (and because of this also $y$) thus has valuation zero. Let $a''$ be the projection of $y'$ on $a$. This element is residually opposite $x''$, but cannot be residually opposite $y''$ as it is only at distance 4 from $y''$. This contradicts the previous case applied to $x''$ and $y''$.  
\end{itemize}
\end{itemize}
\qed

\begin{lemma}\label{lemma:element}
Let $x,y$ be elements of $\Gamma$ such that $[x] \inc_r [y]$. Then there exists $y' \in [y]$ such that $x \inc y'$. 
\end{lemma}
\proof
Let $F$ be the set of all flags containing an element of $[x]$ and one of $[y]$. Let $\{x',y'\}$ be a flag of $F$ such that the sum $d$ of distances of $x'$ and $y'$ to $x$ is minimal. If $d=1$, then $x'=x$ and $x\inc y'$. So we may suppose that $d>1$. 

Suppose that the distance of $x$ to $y'$ is one bigger than the distance from $x$ to $x'$. Let $(x_0=x,x_1,\dots, x_{j-1}=x',x_j=y')$ be the shortest path from $x$ to $y'$ containing $x'$ ($j\leq n$). Let $i$ be the smallest integer such that the subpath $(x_i, \dots, x_{j-1},x_j)$ has valuation zero. We have that $i \geq 1$ (because otherwise it is impossible that $x' \in [x]$) and $i\leq j-1$. Using Lemma~\ref{lemma:pq0} we can extend this subpath to a path $(x_i, \dots, x_{j-1},x_j, x_{j+1},\dots,x_{i+n})$ with valuation zero of length $n$.  Consider the unique path $(x'_i=x_i, x'_{i+1}=x_{i-1}, \dots,x'_{i+n}=x_{i+n})$ from $x_i$ to $x_{i+n}$ containing $x_{i-1}$. Then using (U4), we see that this path has valuation zero. These two paths together form an ordinary $n$-gon $\Omega$, which is simply folded along $x_i$. The previous lemma implies that $x'_{j-1} \in [x]$ and $x'_j \in [y]$. But the sum of distances to $x$ of these two incident elements is strictly less than $d$, contradicting the minimality of $d$.

The case where the distance of $x$ to $x'$ is one bigger than the distance from $x$ to $y$ is proven analogously.
\qed

The diameter of our new geometry $\Gamma_r$ is clearly $n$. In order to prove it is a (weak) generalized $n$-gon we have to show that there is no closed non-stammering path of length less than $2n$. So suppose by way of contradiction that we have such a path $([x_0],[x_1], \dots , [x_{2m}]=[x_0])$ with $2 \leq m < n$. The previous lemma allows us to lift the path into a (not necessarily closed) path $(x'_0 ,x'_1,\dots,x'_{2m})$ such that $[x'_i]=[x_i]$. 

Due to Corollary~\ref{cor:zero} and the fact that the original path was non-stammering, this path has valuation zero. If $2m<n$, we extend this path to a path $(x'_0 ,x'_1,\dots,x'_{2m}, x'_{2m+1}, \dots,x'_{n})$ with valuation zero of length $n$ (this is possible by Lemma~\ref{lemma:pq0}). Whether or not $2m < n$, $x'_n$ is residually opposite $x'_0$, but is not opposite, and thus certainly not residually opposite $x'_{2m}$. Hence we have a contradiction and we have proven the second main result.

\section{Proof of the third main result}
The main idea of the proof is starting from one valuation $u$ on $\Gamma$, to construct more valuations. Each of these valuations will correspond to a point of our $\R$-building. We first cite a lemma from~\cite{kshvm} that we will use in our proof. We use the following notation: the sector-panel with direction $x$ and source $\alpha$ is denoted by $x_\alpha$, the residual distance in the residue of $\beta$ is denoted by $\dd_\beta$, and the length of the intersection of two sector-panels with source $\beta$ and directions $x$ and $y$ as $u_\beta(x,y)$.

\begin{lemma}\label{lemma:progress}
Let $\Lambda$ be an affine apartment system with a generalized polygon $\Lambda_\infty$ at infinity. Let $\alpha$ be a point of $\Lambda$. Let $x,a,b,c$ be elements of $\Lambda_\infty$ such that $a \inc b \inc c$, and $\beta$ a point on $x_\alpha$
with $\dd(\alpha,\beta) =l$. Then there exists $\delta > 0$ such that for any $\beta'$ on $x_\alpha$ with
$\dd(\alpha,\beta') \in [l,l+\delta]$, the following holds~: $$u_{\beta'}(a,c) = u_\beta(a,c) + \epsilon\frac{
\sin(\dd_\beta(b,x_0 ) \pi /n)  }{ \sin(\pi/n)} \dd(\beta,\beta'),$$ where $\epsilon$ is a constant equal to
 $$\begin{cases}
 -1 \mbox{ if }\dd_\beta(a,x ) = \dd_\beta(c,x ) =\dd_\beta(b,x )-1, \\
 1 \mbox{ if } \dd_\beta(a,x ) = \dd_\beta(c,x )=\dd_\beta(b,x )+1, \\
 0 \mbox{ if } \dd_\beta(a,x ) \neq \dd_\beta(c,x ).
\end{cases}$$

\end{lemma}

We now return to our case. Let $(\Gamma,u)$ be a generalized $n$-gon with valuation, $x$ an element of $\Gamma$, and $t \in \R^+$ a positive real number. We want to define a new valuation $u^{V(x,t)}$ with $V(x,t)$ an operator called the \emph{translation operator} ($u^{V(x,t)}$ will be referred to as the \emph{$t$-translation of $u$ towards $x$}, and $u$ is \emph{$t$-translated towards $x$}). 


How do we construct this new valuation? Remember that each element $y$ has a certain residual distance $\dd_r(x,y)$ from $x$ in the residue $\Gamma_r$ defined by $u$. We now `predict' the \emph{translated residual distance} $\dd_r^{x,t}(y)$ from $x$ to $y$ when $t$-translating $u$, as it would be if we were indeed in an affine apartment system (we changed the notation of the residual distance to an unary function to stress the dependability of $x$, and the fact that we will only need distances from $x$). This function defined for $t\in [0,+\infty[$ will be right-continuous and piecewise constant. First thing one needs to assure here is that for two incident elements $y,z$, the translated residual distances $\dd_r^{x,t}(y)$ and $\dd_r^{x,t}(z)$ differ by only one. The definition of this function will be referred to as step (C1), the difference condition as condition (C2).

Because we know how the (translated) residual distances would behave if we were in an affine apartment system, we can use Lemma~\ref{lemma:progress} to predict how the translated individual valuations would behave if we were indeed in an affine apartment system (this is done by a trivial integration of a piecewise constant
function). The set of all these individual valuations allows to construct a new
`valuation' $u^{V(x,t)}$ (we still need to verify this is really a
valuation).
On pages 9-11 of the above mentioned paper~\cite{kshvm} it was shown that the weighted sum of the coefficients of $t$ along the path $(x_0, \dots,x_n)$ depends only on the residual distances of $d_0$ and $d_n$ of $x_0$ and $x_n$ respectively, under the assumption that $d_0=x_0$. The argument in~\cite{kshvm} can be extended to show that this weighted sum depends only on $d_0$ and $d_n$ also when $d_0$ is not zero by applying the same idea as in Case (v) on page 11 of~\cite{kshvm} if $j=1$ is a valley. Because here the predicted individual valuations behave in the same way as they would in the affine apartment system case, this result can be applied here (also using the fact that for two incident elements the residual distances
differ only by one) to guarantee that (U4) will be satisfied by $u^{V(x,t)}$.  The
condition (U2) is trivially satisfied. For more insight in how $u^{V(x,t)}$ is constructed see the example in the section below.


For the other two conditions and positivity of the valuation, we will define and use the $\R$-trees associated to elements of $\Gamma$. Note that an \emph{$\R$-tree} (or simply \emph{tree} when no confusion can arise) is simply an affine apartment system of type $\widetilde{A_1}$, or, equivalently, of dimension one. 

Choose a point $x$ in a given tree. We can define a valuation $v$ acting on the set of pairs $(e,f)$ of ends (parallel classes of sectors) of this tree as the length of the intersection of the 2 half apartments with boundary $x$ and ends $e$ and $f$. The point $x$ will be called the \emph{base point of the valuation}.

One property of $v$ is that if for 3 ends $e,f,g$ the inequality $v(e,f) < v(f,g)$ implies $v(e,g) = v(e,f)$. Now, given any binary function $w$ acting on a set $E$ obeying this property, one can (re)construct a tree (if $w$ is already a valuation of a tree, then we will obtain the same tree) by taking the set $\{(e,t) | e\in E, t\in \R^+\}$ and applying the equivalence relation $$(e,t) \sim (f,s) \Leftrightarrow t=s \mbox{ and }t \leq w(e,f)$$ ($e,f \in E$ and $s,t \in \R^+$). The base point of this tree is the equivalence class  $\{(e,0) | e \in E\}=:x$. The set of ends of this tree is in natural bijective correspondence with $E$ and the valuation in this tree with base point $x$ coincides with $w$. (This construction is a special case of the one of Alperin and Bass in~\cite{Alp-Bas:87}.)

It is easily seen that this property is the same as (U3) when we restrict $u$ to a point row or line pencil. So to each line $L$ or point $p$ of $\Gamma$ we can associate a tree named $T(L)$ or $T(p)$ with a certain base point. The location of this base point will play a major role in the next sections. Other choices of base points yield other valuations of the tree.

We now return to the problem of (U1), (U3) and positivity. Obviously, this will be
solved
if we can show that the change in valuations of elements incident with an
element $y$ of $\Gamma$ is described by changing the base point in the tree
$T(y)$. With an eye on the above lemma, we want to move the base point towards an
end corresponding to an element $a\inc y$ with $d_r^{x,t}(a)=d_r^{x,t}(y)-1$
over a length of $t\sin(d_r^{x,t}(y) \pi /n)  / \sin(\pi/n)$ with $t$ a
certain
translation length such that the translated residual distances of $a$ and
$y$
stay the same. In order that the valuations obtained by this change of
base point correspond to the predictions of the valuations using the above
lemma, we need to verify three things.
\begin{itemize}
\item
If the valuation of the pair consisting of $a$ and another element $b \inc
y$ is
going to decrease (equivalent with saying that $d_r^{x,t}(b)=d_r^{x,t}(y)-1$
and
$d_r^{x,t}(y) \neq n$), then this valuation corresponds with the predicted
valuation using  the displacement of the base point in the tree if the two
half-apartments with ends $a$ and $b$ and source the base point have more 
in common than only the base point, so $u^{V(x,t)}(a,b) > 0$. (We refer to this as condition (C3).)
\item
If the valuation of the pair consisting of $a$ and another element $b \inc
y$ is
going to stay the same (equivalent with saying that
$d_r^{x,t}(b)=d_r^{x,t}(y)+1$), then we have correspondence between the two
predictions if the base point lies in the apartment with ends $a$ and $b$, so
$u^{V(x,t)}(a,b) = 0$. (This will be condition (C4).)
\item
Finally note that if the valuation is going to increase (two elements $b,c
\inc
y$ with $d_r^{x,t}(b)=d_r^{x,t}(c)=d_r^{x,t}(y)+1$), we would need that the
base point lies on the intersection of the apartment with ends $a$ and $b$,
and
the one with ends $a$ and $c$ (thus $u^{V(x,t)}(a,b)  =u^{V(x,t)}(a,c) =
0$).
But this is already covered by (C4), so there is no extra condition needed.
\end{itemize}


In the next part of the proof (after the example), we consider each case seperately. 
\subsection{An example}

We will illustrate with an example how $u^{V(x,t)}$ will be calculated in practice. Suppose we are in the $n=3$ case, and that $x$ is a point. Let us say we have two points $x_1$, $x_2$ different from $x$, and we want to define $u^{V(x,t)}(x_1,x_2)$. (For the (C1) used here, see later on.)

Suppose $u(x,x_i)=t_i$ and suppose $u(x_1,x_2)=t_2$, with $t_1>t_2>0$ (there are other cases, but let's rectrict to this one). The residual distances are all zero between these points. Let $L$ be the line joining $x_1$ and $x_2$. $\epsilon$ in the formula of Lemma~\ref{lemma:progress} is -1. We can take here $\delta=t_2$ (so long, the residual distances to x do not change according to (C1)),  and we obtain
\begin{equation}
u^{V(x,t)}(x_1,x_2)=t_2-t \mbox{ for } t\leq t_2
\end{equation}
From then on, $\epsilon$ becomes zero until $t_1$, since the residual distances to $x$ from $x_1$ differs from that to $x_2$; to $x_2$ it becomes 2 and to $x_1$ it is 0. Hence
\begin{equation}
u^{V(x,t)}(x_1,x_2)=0 \mbox{ for } t_2  < t \leq t_1
\end{equation}
Note that, up to now, the residual distance from $x$ to $L$ was always 1, hence the quotient of the $\sin$'s has always been 1. This is going to change now.

After that, $\epsilon$ becomes 1, and the quotient of the $\sin$'s is still 1, but only until $\tau(x,L)$ according to (C1), which is by definition bigger than $t_1$. Hence
\begin{equation}
u^{V(x,t)}(x_1,x_2)=t-t_1 \mbox{ for }  t_1 < t \leq \tau(x,L)
\end{equation}
At $t=\tau(x,L)$, the $\sin$ of $\dd(x,L)\pi/3$ becomes 0, and so the valuation becomes constant again:
\begin{equation}
u^{V(x,t)}(x_1,x_2)=\tau(x,L)-t_1 \mbox{ for }  \tau(x,L) < t. 
\end{equation}

\subsection{$n=3$}
We check (C1), (C2), (C3) and (C4).
\subsubsection{(C1)}
\begin{itemize}
\item If $\dd(x,y) = 0$, then $\dd_r^{t,x}(y)=0$ for $t \in [0, +\infty[$.
\item If $\dd(x,y) = 1$, then $\dd_r^{t,x}(y)=1$ for $t \in [0, +\infty[$.
\item If $\dd(x,y) = 2$, then 
\begin{itemize}
\item $\dd_r^{t,x}(y)=0$ for $t \in [0,  u(x,y)[$,
\item $\dd_r^{t,x}(y)=2$ for $t \in [u(x,y), +\infty[$.
\end{itemize}
\item If $\dd(x,y) = 3$, then
\begin{itemize}
\item $\dd_r^{t,x}(y)=1$ for $t \in [0,  \tau(x,y)[$, 
\item $\dd_r^{t,x}(y)=3$ for $t \in [\tau(x,y), +\infty[$.
\end{itemize}
\end{itemize}

\subsubsection{(C2)}
Let $y$ and $z$ be a pair of incident elements. Without loss of generality one can suppose that $\dd(x,y)+1 = \dd(x,z)$. The only not completely trivial cases is where $\dd(x,y) =2$ and $\dd_r^{t,x}(y)=0$. This happens when $t \in [0,  u(x,y)[$, so also $t < \tau(x,z) = u(x,y)+u(y,z)$, and thus $\dd_r^{t,x}(z)=1$. We conclude that (C2) is satisfied.

\subsubsection{(C3)}
Let again $y$ be an element, with $a,b$ two elements incident with $y$, such that $\dd_r^{x,t}(a) + 1 = \dd_r^{x,t}(b) + 1 = \dd_r^{x,t}(y)$. The only cases for which we need to verify that $u^{V(x,t)}(a,b) > 0$ are $\dd_r^{x,t}(y) = 1$ or $2$.
\begin{itemize}
\item If $\dd(x,y) = 1$, then $\dd_r^{x,t}(a) + 1 = \dd_r^{x,t}(b) + 1 = \dd_r^{t,x}(y)=1$. One can choose $a=x$, then $\dd(x,b) = 2$, so in this case $t \in [0,  u(x,b)[ $. The following now holds: $u^{V(x,t)}(a,b) = u(x,b) - t > 0$. 
\item If $\dd(x,y) = 2$, then $\dd_r^{t,x}(y)=2$ for $t \in [u(x,y), +\infty[$. Assume that $a=xy$ and $\dd(x,b)=3$. This yields that $t \in [u(x,y), \tau(x,b)[ = [u(x,y),u(x,y) + u(a,b)[$. One checks that $u^{V(x,t)}(a,b) = u(a,b) -t + u(x,y) >0$, so (C3) holds here.
\item If $\dd(x,y) = 3$, then $\dd_r^{t,x}(y)=1$ for $t \in [0,  \tau(x,y)[$. This case is similar to the case $\dd(x,y) = 1$, but now using Lemma~\ref{lemma:path} instead of (U3).
\end{itemize}  
\subsubsection{(C4)}
Let $y$ be an element, with $a,b$ two elements incident with $y$, such that $\dd_r^{x,t}(a) + 1 = \dd_r^{x,t}(b) - 1 = \dd_r^{x,t}(y)$. We only need to verify that $u^{V(x,t)}(a,b) = 0$ is when $\dd_r^{x,t'}(b) <\dd_r^{x,t}(b)$ for $t'<t$. 
\begin{itemize}
\item If $\dd(x,y) = 1$, we again choose $x$ to play the role of $a$. It is clear that the conditions then tell that $t=u(x,b)$, and $u^{V(x,t)}(x,b) = u(x,b)-t =0$.
\item If $\dd(x,y) = 2$, then $\dd_r^{t,x}(y)=2$ for $t \in [u(x,y), +\infty[$. We choose $a$ to be the element $xy$. The element $b$ thus lies at distance 3 from $x$, and $t=\tau(x,b)$. Similarly to the $(C3)$ case one checks that $u^{V(x,t)}(a,b) = u(a,b) -t + u(x,y) =0$.
\item If $\dd(x,y) = 3$, then $\dd_r^{t,x}(y)=1$ for $t \in [0,  \tau(x,y)[$. This case is similar to the case $\dd(x,y) = 1$, but now using Lemma~\ref{lemma:path} instead of (U3).

\end{itemize}
\subsection{$n=4$}
Before we check the conditions, we state some useful lemmas.
\begin{lemma}\label{lemma:oddquad}
It is impossible to have an ordinary quadrangle $\Omega$ containing exactly two sides with non-zero valuations, such that opposite elements have the same valuation, but each two corners of a side have different valuations.
\end{lemma}
\proof
Suppose that such a quadrangle $\Omega$ does exist. Then let $p,q$ be corners of $\Omega$ such that $u(p,q)>0$, and such that the valuation in $p$ is bigger than the one in $q$. There exists an $r \inc pq$ such that $u(p,r)=u(q,r)=0$ (by Lemma~\ref{lemma:pq0} and (U3)). Let $\Omega'$, $\Omega''$ be the ordinary quadrangles sharing a path of length 4 with $\Omega$ and containing $r,p$ and $r,q$, respectively. Denote the element opposite $pq$ in $\Omega$ by $s$. Let $p'$, $q'$ and $r'$ be the projections of, respectively, $p$, $q$ and $r$ on $s$. Because the valuation in $p$ is bigger than the one in $q$, (U4) applied in both $\Omega'$ and $\Omega''$ yields $u(r',q') < u(r',p')$ (because these are the only two other terms in applying (U4) differing in both quadrangles), so $u(r',q')=u(p',q') > 0$ by (U3).

The valuations of the elements $r$ and $r'$ in $\Omega'$ cannot be equal because the valuation of $q$ in $\Omega'$ is strictly smaller than the valuation of $q'$ in $\Omega'$. So the two corners with smallest valuation in $\Omega'$ --- guaranteed by (the dual of) Lemma~\ref{lemma:path} ---  have to be in the corners $q$ and $r'$. Applying (U4) we obtain $u(q,q')+ \sqrt{2} u(qq',qr) + u(q,r) = u(q',r')+\sqrt{2} u(r'q',r'r) +u(r,r')$, which implies that $u(q',r')=0$, a contradiction.
\qed

\begin{lemma}\label{lemma:samedist}
Let $a,b$ be two opposite elements. Then there exist two paths $(a,x_1,x_2,x_3,b)$ and $(a,y_1,y_2,y_3,b)$ from $a$ to $b$ such that $u(a,x_2) = u(x_2,b)$, $u(a,y_2) = u(y_2,b)$ and $u(x_1,y_1) =0$, if and only if for each path $(a,z_1,z_2,z_3,b)$ the equality $u(a,z_2) = u(z_2,b)$ holds.
\end{lemma}
\proof
The implication from right to left is trivial by (U1). So suppose the left part of the statement is satisfied.

First remark that (U4) tells us that $u(x_3,y_3)=0$, so the situation is symmetric in $a$ and $b$. Suppose that $u(a,z_2) < u(z_2,b)$, then without loss of generality we may assume that $u(x_1,z_1) = 0$ (by (U3)). But then $u(x_2,a)+\sqrt{2} u(x_1,z_1) + u(a,z_2) < u(x_2,b)+\sqrt{2} u(x_3,z_3) + u(b,z_2)$, which contradicts (U4).
\qed

If for two opposite elements $a$ and $b$ the situation of the above lemma holds, then we say that those two points are \emph{equidistant}.

\begin{lemma}\label{lemma:choice}
If two opposite points $x,y$ are not equidistant, then there exists a path $(x,a,b,c,y)$ from $x$ to $y$, such that $u(x,b) \geq u(b,y)$ and $u(a,c)=0$.
\end{lemma}
\proof
First note that, if for all paths $(x,a',b',c',y)$ from $x$ to $y$ it would happen that $u(x,b') \leq u(b',y)$,  then condition (U4) or Lemma~\ref{lemma:samedist} is violated in a quadrangle defined by two paths $(x,a',b',c',y)$ and $(x,a'',b'',c'',y)$, where $a'$ and $a''$ are chosen so that $u(a',a'')=0$ (which is possible due to (U1)). 

So we know the existence of a path $(x,a',b',c',y)$ with $u(x,b') > u(b',y')$. If $u(a',c')=0$,  then we are finished, so assume this is not the case. Using Lemma~\ref{lemma:pq0}, we can find $a'' \inc x$ with $u(a',a'')=0$. Let $(x,a'',b'',c'',y)$ be the unique shortest path from $x$ to $y$ containing $a''$. Lemma~\ref{lemma:path} tells us that either $u(c',c'')=0$ or $u(a'',c'')=0$. If we are in the first case, then applying Lemma~\ref{lemma:path} again on the other type of elements in the ordinary quadrangle leads to a contradiction with Lemma~\ref{lemma:oddquad}. So $u(a'',c'')=0$. Using (U4) one sees that $(x,a'',b'',c'',y)$ is a path with the desired properties.
\qed

We are now ready to check (C1), (C2), (C3) and (C4).
\subsubsection{(C1)}
\begin{itemize}
\item If $\dd(x,y) = 0$, then $\dd_r^{t,x}(y)=0$ for $t \in [0, +\infty[$.
\item If $\dd(x,y) = 1$, then $\dd_r^{t,x}(y)=1$ for $t \in [0, +\infty[$.
\item If $\dd(x,y) = 2$, then 
\begin{itemize}
\item $\dd_r^{t,x}(y)=0$ for $t \in [0,  u(x,y)[$,
\item $\dd_r^{t,x}(y)=2$ for $t \in [u(x,y), +\infty[$.
\end{itemize}
\item If $\dd(x,y) = 3$, with $x\inc a \inc b \inc y$ then
\begin{itemize}
\item $\dd_r^{t,x}(y)=1$ for $t \in [0,  u(x,b)  + u(a,y)/\sqrt{2}[$, 
\item $\dd_r^{t,x}(y)=3$ for $t \in [u(x,b)  + u(a,y)/\sqrt{2}, +\infty[$.
\end{itemize}
\item If $\dd(x,y) = 4$, then in the case that there exist $a,b$ and $c$ such that $x\inc a \inc b \inc c \inc y$, with $u(x,b) \neq u(b,y)$, let $k(x,y)$ be the minimum of both (this is independent of $a,b$ and $c$ due to Lemma~\ref{lemma:path}). In the case that $x$ and $y$ are equidistant, we define $k(x,y)$ to be equal to $\tau(x,y) / 2$.
\begin{itemize}
\item $\dd_r^{t,x}(y)=0$ for $t \in [0,k(x,y) [$, 
\item $\dd_r^{t,x}(y)=2$ for $t \in [k(x,y),\tau(x,y)-k(x,y)[$.
\item $\dd_r^{t,x}(y)=4$ for $t \in [\tau(x,y) - k(x,y),+\infty[$.
\end{itemize}
\end{itemize}

\subsubsection{(C2)}
Let $y,z$ be a pair of incident elements. Without loss of generality one can suppose that $\dd(x,y)+1 = \dd(x,z)$. There are three nontrivial cases.

\begin{itemize}
\item $\dd(x,y) = 2$, with $\dd_r^{t,x}(y)=0$, and $\dd_r^{t,x}(z)=3$. This yields $t \in [0,u(x,y)[ \cap [u(x,y)+u(xz,z)/\sqrt{2}, +\infty[$. The last intersection is clearly empty and thus this case cannot occur.

\item $\dd(x,y) = 3$, with $\dd_r^{t,x}(y)=1$ and $\dd_r^{t,x}(z)=4$. Let $x\inc a\inc b\inc y$. This situation occurs when $t \in  [0,  u(x,b)  + u(a,y)/\sqrt{2}[ \cap [\tau(x,z) - k(x,z),+\infty[$. As $k(x,z) \leq \min(u(x,b),u(b,z))+ u(a,y) )/\sqrt{2} $, the range for $t$ is empty, thus this case cannot occur either. 

\item $\dd(x,y) = 3$, with $\dd_r^{t,x}(y)=3$ and $\dd_r^{t,x}(z)=0$. Let $x\inc a\inc b\inc y$. This happens when $t \in  [0,  k(x,z)[ \cap [u(x,b)  + u(a,y)/\sqrt{2}, +\infty[$.  Again the bound $k(x,z) \leq \min(u(x,b),u(b,z))+ u(a,y) )/\sqrt{2}$ leads to a contradiction.

\end{itemize}

\subsubsection{(C3)}
Let again $y$ be an element, with $a,b$ two elements incident with $y$, such that $\dd_r^{x,t}(a) + 1 = \dd_r^{x,t}(b) + 1 = \dd_r^{x,t}(y)$.
\begin{itemize}
\item If $\dd(x,y) = 1$, then $\dd_r^{t,x}(y)=1$ for $t \in [0, +\infty[$. Let $a$ be the element $x$.  then $\dd(x,b) = 2$, so in this case $t \in [0,  u(x,b)[ $. The following now holds: $u^{V(x,t)}(a,b) = u(x,b) - t > 0$. 
\item If $\dd(x,y) = 2$, then $\dd_r^{t,x}(y)=2$ for $t \in [u(x,y), +\infty[$. We may assume that $a=xy$ and $\dd(x,b)=3$. This yields that $t \in [u(x,y), [u(x,y),u(x,y) +  u(a,b)/\sqrt{2}[$. One checks that $u^{V(x,t)}(a,b) = u(a,b) -\sqrt{2} (t - u(x,y)) >0$, so (C3) holds here.

\item If $\dd(x,y) = 3$, with $x\inc p \inc q \inc y$ then
\begin{itemize}
\item $\dd_r^{t,x}(y)=1$ for $t \in [0,  u(x,q)  + u(p,y)/\sqrt{2}[$. We distinguish two subcases.
\begin{itemize}
\item If $u(x,q) > t$, then we choose $a=q$. The element $b$ is then at distance 4 from $x$, with $\dd_r^{t,x}=0$, hence $t \in [0, k(x,b)[$. If $u(q,b)\leq t$, then $u(q,b) = k(x,b) \leq t$ which is impossible (remember $u(x,q) > t$). As $u^{V(x,t)}(q,b)=u(q,b) -t$, condition (C3) is satisfied here.
\item The other subcase is where $u(x,q) \leq t$. Note that $\dd_r^{t,x}=2$, thus $\dd(x,b) =4$. Since $u(x,q) \leq t$ and $t < k(x,b)$, we have $u(q,b)=u(x,q)$ and $u(p,y) >0$. We construct $a$ as follows: let $r$ be an element incident with $x$ such that $u(p,r)=0$ and let $s$ be an element incident with $r$ such that $u(x,s)=0$. The element $a$ is the projection of $s$ on $y$. Let $c$ be the projection of $b$ on $r$. Lemmas~\ref{lemma:path} and~\ref{lemma:oddquad} yield $u(a,s)=u(y,as)=0$, $u(r,as) =\tau(x,a)/ \sqrt{2} $, $a$ and $x$ are equidistant (by Lemma~\ref{lemma:samedist}), and $\dd_r^{t,x}(x,a)=0$. As $u^{V(x,t)}(a,b) = u(a,b)-t$, we have to prove that $u(a,b) \geq k(x,b)$ in order to prove (C3).  

Let $\Omega$ be the unique quadrangle containing $b$, $y$, $s$ and $r$. If $b$ and $x$ are equidistant, then the valuation of $b$ in $\Omega$ is zero, and (U4) implies $u(a,b) \geq u(r,as)/\sqrt{2}= k(x,b)$. Finally suppose that $b$ and $x$ are not equidistant, then Lemma~\ref{lemma:samedist} implies $u(x,s) \neq u(s,c)$, and thus $u(x,s),u(s,c) \geq k(x,b)$ (by definition of $k(x,b)$). Applying (U4) in $\Omega$ tells us now that $u(a,b) \geq u(s,c) \geq k(x,b)$, which we needed to show.

\end{itemize}
\item $\dd_r^{t,x}(y)=3$ for $t \in [u(x,q)  + u(p,y)/\sqrt{2}, +\infty[$. Let $a$ be $q$ in this case. The element $b$ will thus be at distance $4$, while $\dd_r^{x,t}(b)=2$. So $t \in [k(x,b),\tau(x,b)-k(x,b)[$, which also means that $b$ and $x$ are not equidistant. Careful analysis reveals that $u^{V(x,t)}(a,b) = \tau(x,b) -k(x,b) -t$, which is strictly larger than zero because $\dd_r^{x,t}(b)=2$ implies that $t \in [k(x,b), \tau(x,b) -k(x,b)[$.

\end{itemize}
\item If $\dd(x,y) = 4$, then $\dd_r^{t,x}(y)=2$ for $t \in [k(x,y),\tau(x,y)-k(x,y)[$. Notice that $x$ and $y$ are not equidistant. Let $(x,p,q,a,y)$ be a path as constructed in Lemma~\ref{lemma:choice}. This fixes our choice of $a$. Let $(x,r,s,b,y)$ be the unique path from $x$ to $y$  containing $b$. One checks that $u^{V(x,t)}(a,b) = u(a,b) - \sqrt{2}(t-k(x,y))= u(a,b) - \sqrt{2}(t-u(y,q))$. The value of $t$ is strictly smaller than $u(x,s)+u(s,b)/\sqrt{2}$ (because $\dd_r^{t,x}(b)=1$). All we have to check is that $u^{V(x,t)}(a,b) \geq 0$ when $t =u(x,s)+u(s,b)/\sqrt{2}$. Using (U4), one proves that $u^{V(x,t)}(a,b) = u(p,r) \geq 0$ for this value of $t$. 
\end{itemize}
This concludes the proof of (C3) in this case.

\subsubsection{(C4)}
The condition (C4) can be proved analogously to the proof of (C3) in this case.

\subsection{$n=6$ and the valuation is discrete}
Here the discreteness allows us to define the translations in a much easier way using recursion. We start with a valuation $u$ where the valuations of one type of elements are integer multiples of 3, while valuations of the other type are integer multiples of $\sqrt{3}$ (with proper rescaling, this is a consequence of the discreteness, see~\cite{kshvm}). The valuation $u$ also defines a residual distance $\dd_r$. We use this as a constant translated residual distance $\dd_r^{x,t}$ with $t \in [0,1[$ or $[0,\sqrt{3}/2[$ depending on the type of $x$ (notice that (C1) and (C2) are fulfilled by this). The condition (C4) is satisfied because it is satisfied for $t=0$, and that the valuations in question stay zero. The discreteness makes it so that because (C3) is satisfied for $t=0$, it will also be satisfied for $t$ in the ranges above (because the range is small enough so the valuation in question cannot decrease to zero).

Let's clarify this with an example first. Suppose that $x$ is an element such the valuations of that type of element are integer multiples of $\sqrt{3}$, and let $k \in [0,\sqrt{3}/2]$. Applying what is said above the displacement of the base point of the trees associated with element $y$ with residual distance $\dd_r(x,y)$ to yield the valuation $u^{V(x,k)}$ will be as given in the following table, all displacements are towards an element which is in the residue closest to $x$:

\begin{table}[htdp]
\begin{center}\begin{tabular}{c|c|c|c|c|c|c|c}$\dd_r(x,y)$ & 0 & 1 & 2 & 3 & 4 & 5 & 6 \\\hline Displacement of basepoint & none & $k$ & $\sqrt{3}k$ & $2k$ & $\sqrt{3}k$ & $k$ & none\end{tabular} 
\end{center}
\label{defaulttable}
\end{table}

Note that $k$ is small enough so that the displacements don't make the base points reach branching points of the trees, except for the maximal value $k=\sqrt{3}/2$ and $\dd_r(x,y)=3$. Branching points may not be crossed, because for (C3) valuations may not decrease to zero (which is what happens at branching points), except for the final point (for a $k$-translation, (C3) needs only to be checked for values $t$ in $[0,k[$).

We can repeat the same procedure on the new valuations we obtain but with one major caveat: the valuations are no nice integer multiples anymore (because we can $k$-translate with $k$ a real number in $ [0,1]$ or $[0,\sqrt{3}/2]$ depending on the situation). However, this does not pose an unsolvable problem. Let $W$ be a Coxeter group of type $\widetilde{\mathsf{G}}_2$ acting naturally on an Euclidean affine plane $\A$. Take a special vertex $s$. Notice that, with proper rescaling, the distances from $s$ to all the walls of a parallel class of walls is exactly the image set of the valuations $u$ of the elements incident with a certain type of elements. Let $s'$ be a point of the plane $\A$ at distance $k$ from $s$, on the same wall (with type the element we have translated to) as $s$.  Due to Lemma~\ref{lemma:progress} (or by looking at the example above), we can again identify distances from $s'$ to all the walls of a parallel classes to image sets of valuations $u^{V(x,k)}$ of certain elements as above. (We can no longer identify with a type of elements, there will be more classes of elements, due to the residue corresponding with $u^{V(x,k)}$ being a weak generalized hexagon.)

We can now $l$-translate $u^{V(x,k)}$ to an element $y$ in the same way as above, with $l$ small enough that we don't `cross' any walls with the corresponding displacement of the point in the plane. The displacement will now along the line at angle $d\pi/n$ with the line through $s$ and $s'$, with $d$ the distance in the residue of $u^{V(x,k)}$ from $x$ to $y$. One cannot cross the wall because we will have moved some base points of trees to branching points. Note however that `arriving' at a wall is allowed, so one can get across that wall with the next translation).  

This procedure allow us to repeat the construction, obtaining all subsequent translations of $u$ we want.

We again clarify further with an example. Suppose $x$ as in the above example and let $t$ be $\sqrt{3}/3$. Now suppose that $y$ is an element which is at distance 2 from $x$ in the residue of $u^{V(x,k)}$. With the above procedure it follows that we $l$-translate to $y$ with $l \in [0,\sqrt{3}/3]$ (when $l=\sqrt{3}/3$, we arrive again in a special point of $\A$). Again we could make a table and confirm that indeed that the base points reach branching points of the tree except for the maximal value $l=\sqrt{3}/3$.

\subsection{What about $n=5$ and the non-discrete case for $n=6$?}
One could use similar techniques as for the cases $n=3$ and $n=4$ to investigate these cases. The things one would need to prove are mostly quantitative versions of the qualitative lemmas of the proof of the first main result. However extending the, already extensive, complexity of the case studies $n=3$ and $n=4$ to these higher cases, would probably require a comparable number of pages than the current version of the paper. For this reason we choose to restrict ourselves to the already handled cases.

\subsection{Some first observations}
Now that we defined additional valuations, we need to show that they form the point set of an $\R$-building. We need some properties to do so. 

\begin{lemma}\label{lemma:resequiv}
The residual distance of $x$ and $y$ in the residue of $u^{V(x,t)}$ equals $\dd_r^{t,x}(y)$.
\end{lemma}
\proof 
This follows from the way we defined (C1) for $n=3$ and $n=4$, and from the construction for the discrete case when $n=6$.
\qed

\begin{lemma}
If $\dd_r^{x,t}(y)=n$, then $\dd_r^{x,t'}(y)=n$ for every $t' \geq t$.  
\end{lemma}
\proof
This only case this isn't directly clear is $n=6$. Applying the previous lemma we see that in the residue of $u^{V(x,t)}$ the elements $x$ and $y$ are residually opposite and that each shortest path between both has valuation zero. Because of the way we defined $u^{V(x,t')}$, it follows that the path also has valuation zero for $u^{V(x,t')}$. This proves the lemma.
\qed  

\begin{cor}\label{cor:inc}
When translating towards $x$, the residual distance $\dd_r^{x,t}(y)$ only increases, up to the point that $\dd_r^{x,t}(y)=\dd(x,y)$.
\end{cor}
\proof
Again we only need to prove this when $n=6$. Because of the previous lemma and the fact that the residue is a weak generalized $n$-gon where each element is incident with at least 2 elements, we see that $\dd_r^{x,t}(y)$ only increases. It increases to $\dd(x,y)$ because if for an arbitrary element $z$ we have $\dd_r^{x,t}(z)=\dd(x,z)<n$, then for an element $a \inc z$ there exists $t'\geq t$ such that $\dd_r^{x,t}(a)=\dd(x,a)$ (this is due to the displacement of the basepoint of the tree associated to $z$, which happens at a constant rate towards the projection of $x$ on $z$). Repeating this argument implies that $\dd_r^{x,t}(y)$ will eventually become $\dd(x,y)$.
\qed




\subsection{Structural properties of the set of translated valuations}
Let $\Lambda(u)$ be the set of all valuations obtained by translating $u$ a finite number of times.


\begin{lemma}
If we know the values of a valuation $v$ on the pairs of elements incident with an element $x$, and we know that an element $y$ is residually opposite $x$, then we know the values of $v$ on the pairs of elements incident with $y$.
\end{lemma}
\proof
Let $a,b \inc y$, then (U4) in an $n$-gon containing $a,b,x$ and $y$ tells us that $v(a,b)=v(a',b')$ where $a'$ and $b'$ are the projections on $x$ of $a$ and $b$, respectively.
\qed

\begin{lemma} \label{lemma:known}
Let $\Omega$ be an $n$-gon in $\Gamma$, non-folded for a valuation $v \in \Lambda(u)$, such that all values of $v$ in the line pencils of the corners and points on the sides of $\Omega$ are known, then the values of $v$ are known entirely.
\end{lemma}
\proof
Let $x$ be an element of $\Gamma$. Let $y$ be an element of $\Omega$ with minimal distance $k$ to $x$. Notice that $k<n$. If $k=0$, then we know the valuations of pairs of elements incident with $x$, so suppose $k>0$. Let $z$ be the projection of $x$ on $y$. Then there are 2 ordinary $n$-gons containing $z$ and sharing a path of length $n$ with $\Omega$. By applying (U3), (U4) at least one of these two $n$-gons is non-folded for the valuation $v$. Let $\Omega'$ be such an $n$-gon.  The valuations in the line pencils of the corners and points on the sides of $\Omega'$ are known because of the previous lemma. The minimal distance from $x$ to an element of $\Omega'$ is now strictly less than $k$. So by repeating the above argument one sees that one knows the value of $v$ everywhere.
\qed

\begin{cor}\label{cor:equiv}
If $\dd_r^{t',x}(y)=0$ for all $t' \in [0,t[$, then $u^{V(x,t)}=u^{V(y,t)}$.
\end{cor} 
\proof
If $n=6$, then this follows from the `discrete' construction.

In the other cases, let $\Omega$ be a non-folded $n$-gon (for $u$) containing $x$. If we can prove that for each element $z$ in $\Omega$ the relation $\dd_r^{t',x}(z)=\dd_r^{t',y}(z)$ holds for all $t' \in [0,t[$, then the displacements of the base points in the trees corresponding to the elements of $\Omega$ are the same, so by the previous lemma also $u^{V(x,t)}=u^{V(y,t)}$. Moreover, it suffices to prove this for $z$ equal to $x$ and equal to the element opposite $x$ in $\Omega$ because of (C2). 

If $z=x$, then note that, due to the symmetry of the conditions in (C1),  $\dd_r^{t',y}(x)=0$ is equivalent with $\dd_r^{t',x}(y)=0$ for all $t'\in \mathbb{R}^+$, so also for $t' \in [0,t[$. So the result follows from the assumption. 

If $z$ is opposite $x$ in $\Omega$, note that due to the residual equivalency of $x$ and $y$ (by Lemma~\ref{lemma:resequiv}), we have that $\tau(x,z)= \tau(y,z)=0$, and thus $\dd_r^{t',x}(z)=\dd_r^{t',y}(z) =n$ for all $t'\in \mathbb{R}^+$.
\qed

\begin{remark}\rm
It should also be noted that at this point one can prove that the group of projectivities of a line $L$ preserves the tree structure associated with $L$. This allows for a characterization due to Jacques Tits in the case $n=3$, which was formulated without proof in \cite{tits84}.
\end{remark}

\subsection{Apartments}
An apartment in our $\R$-building will consist of all valuations in $\Lambda(u)$ for which a given ordinary $n$-gon is non-folded. Here, we investigate which valuations keep a given ordinary $n$-gon non-folded.  This will give us later on the affine structure of the apartments. 

Let $u$ be a valuation, and let $\Omega$ be a non-folded $n$-gon in $\Gamma$ containing an element $x$. Note that due to (U4) and multiple use of Lemma~\ref{lemma:pq0} each flag can be embedded in such a non-folded $n$-gon, so results obtained here for single points or flags of $\Omega$ are true for all points or flags.

Using the definition of $t$-translation one easily obtains that a translation $V(x,t)$ moves the basepoint of the tree corresponding to an element $y$ of $\Omega$ along the apartment of that tree with ends the two elements of $\Omega$ incident with $y$. The new basepoint lies at length $t \sin(\dd(x,y) \pi/n)/ \sin(\pi /n)$ towards the projection of $x$ on $y$ (note that in the cases that this projection is not defined, the length will be zero). 

Consider the affine real two-dimensional space $\A$. One can think of this as a (degenerate) affine apartment system with an ordinary $n$-gon at infinity. Identify this $n$-gon with $\Omega$ and let $\alpha$ be a point of $\A$. Now consider the point at distance $t$ on the sectorpanel with source $\alpha$ and direction $x$. An important observation that for an element $y$ of $\Omega$ at infinity, the distance component from the original to the new point perpendicular to the direction to $y$ is $t \sin(\dd(x,y) \pi/n)/ \sin(\pi /n)$ exactly the same as above. 

Note also that  $\Omega$ is non-folded for the valuation $u^{V(x,t)}$, and that the displacement of the base points in the aforementioned trees describe $u^{V(x,t)}$ completely when $u$ is known, due to Lemma~\ref{lemma:known}. So we can identify the points of $\A$ with the valuations obtained by translating $u$ to elements of a certain non-folded $n$-gon for $u$. This spawns a few direct consequences.

\begin{cor} \label{cor:apart}
Let $x$ be an element of $\Gamma$ and let $t$ and $s$ be non-negative real numbers. Then
\begin{itemize}
\item $u^{V(x,t) V(x,s)} = u^{V(x,t+s)}$ (local additivity).
\item $u^{V(x,t) V(y,s)} = u^{V(y,s) V(x,t)}$ if $x \inc y$ (local commutativity).
\item $u^{V(x,t) V(y,t)} = u$ if $\tau_u(x,y)=0$ (reversibility).
\item If a path $(x_0,x_1,\dots,x_i)$ (with $i \leq n$) has valuation zero for some valuation $u$, and suppose that $v$ is a valuation obtained from $u$ by subsequently translating towards the respective elements of the path. Then there exists a $j \in \{1, \dots, i\}$ and $t',s' \in \R^+$ such that $v=u^{V(x_{j-1},t') V(x_j,s')}$. In addition, the total sum of lengths of all the translations does not increase. 
\end{itemize}
\end{cor}
Note that the reversibility statement also implies that, if $v \in \Lambda(u)$, then $\Lambda(v)=\Lambda(u)$.

\subsection{Convexity}
The next thing to investigate is how an ordinary $n$-gon $\Omega$ behaves with respect to translations towards elements outside $\Omega$. This will allow us to prove the (convexity) condition (A2) later on.
\begin{lemma}\label{lemma:conv1}
Let $\Omega$ be an ordinary $n$-gon and $x$ an element not residually equivalent to any of the elements of $\Omega$. Then $\Omega$ cannot be a non-folded $n$-gon for $u^{V(x,t)}$ with $t>0$.
\end{lemma}
\proof
Consider the closed path $(x_0, ... , x_{2n}= x_0)$ that $\Omega$ forms. There is an $i \in \{1,\dots,2n\}$ such that the residual distances from $x$ to $x_{i-1}$ and $x_{i+1}$, are both larger than the residual distance from $x$ to $x_{i}$. We excluded that $x_i$ is residually equivalent $x$, so the right derivative (to $t$) of the valuation
$u^{V(x,t)}(x_{i-1},x_{i+1})$ is positive in a certain interval (for $t$) containing 0 where
the residual distances to $x$ in the path are constant. This implies that $\Omega$ is
not non-folded for $t$ in this interval different from zero. We also know that we
can partition $[0,+\infty[$ in a finite set of intervals with constant residual
distances to x in the path, so repeating the above argument proves the lemma. 
\qed

\begin{lemma} \label{lemma:convex}
Let $\{p,L\}$ be a flag in $\Gamma$, let $l,m$ be positive real numbers, and let $\Omega$ be a non-folded $n$-gon. Then, if $\Omega$ is non-folded for the valuation $u^{V(p,l) V(L,m)}$, it is also non-folded for the valuations $u^{V(p,l') V(L,m')}$, for all $l' \in [0,l]$ and $m' \in [0,m]$. Moreover, there is a point $p'$ and line $L'$ in $\Omega$ such that $u^{V(p,l') V(L,m')} = u^{V(p',l') V(L',m')}$ for all $l' \in [0,l]$ and $m' \in [0,m]$.
\end{lemma}
\proof
First for the part till `moreover': using Corollary~\ref{cor:inc} it follows that if we are translating to a certain flag $\{p, L\}$, we can first `use up' that much of the translations to $p$ and $L$ (note that these commute) such that we only
end up with valuations to elements not residually equivalent to an element of
the ordinary $n$-gon. If we now translate further than this, the apartment loses
his non-foldedness and never regains it, due to Lemma~\ref{lemma:conv1}. So if for $u^{V(p,l) V(L,m)}$ the $n$-gon $\Omega$ is still non-folded, it has to be that $p$ and $L$ stay residually equivalent
to elements of the $n$-gon for the whole translation. So if we translate `less' ($u^{V(p,l') V(L,m')}$ with $l' \in [0,l]$ and $m' \in [0,m]$), $\Omega$ will still be non-folded.

The second part now follows from Lemma~\ref{lemma:conv1} and Corollary~\ref{cor:equiv} (the elements $p$ and $L$ stay residually equivalent to the same pair of incident elements of the $n$-gon for
the whole translation because of Corollary~\ref{cor:inc}).
\qed

\subsection{Existence of apartments containing two valuations}

\begin{lemma} \label{lemma:2ok}
Let $u$ be a valuation, and $v,w  \in \Lambda(u)$. Then there exists a point $p$ and line $L\inc p$ in $\Gamma$, and non-negative real numbers $k$ and $l$ such that $w=v^{V(p,k) V(L,l)}$.
\end{lemma}
\proof
First remark that $w \in \Lambda(u) = \Lambda(v)$. So $w$ can be obtained from $v$ with a series of $i$ translations. We prove with induction that this series of translations can be reduced into the desired form.

If $i \leq  1$ this is trivial. If $i >1$ we can reduce the last $i-1$ translations into the desired form, so we have that $w=v^{V(x,k) V(y,l) V(z,m)}$ with $y \inc z$ and $k,l,m \in \R^+$ (note that the last two translations commute). 

We now start a second induction on $j=\max(\dd(x,y),\dd(x,z))$. If this is $1$, then we are done because of Corollary~\ref{cor:apart}. So suppose that $j>1$, and that we can reduce to the desired form if the maximum is strictly less than $j$. Without loss of generality, assume the maximum in the definition is reached for $\dd(x,z)$. Let $t$ be the smallest real positive number such that the residual distance between $x$ and $z$  in $v^{V(x,t)}$ equals the actual distance in $\Gamma$. There exists an element $x'$ such that $\dd(x',z) < \dd(x,z)$ and $x'$ is residually equivalent with $x$ for $v^{V(x,t')}$, with $t'<t$ (the existence of such an $x'$ will be clarified below). 

If $k \leq t$, then $w=v^{V(x,k) V(y,l) V(z,m)}=v^{V(x',k) V(y,l) V(z,m)}$, and so we are done in this case by the second induction hypothesis. If $k > t$, then $$w=v^{V(x,k) V(y,l) V(z,m)}=(v^{V(x,t)})^{V(x,k-t) V(y,l) V(z,m)}.$$ By the definition of $t$, there exists a non-folded $n$-gon for the valuation $v^{V(x,t)}$ containing $x,y$ and $z$. This implies that the last three translations can be reduced in the desired form of two translations towards two incident elements in the path from $x$ to $z$ (by Corollary~\ref{cor:apart}). If both of these translations are not towards $z$, then we are done due to the second induction hypothesis. If this is not the case then  $w=(v^{V(x,t)})^ {V(y,l') V(z,m')}=(v^{V(x',t)})^ {V(y,l') V(z,m')}$ for certain $l'$ and $m'$, which is again reducable due to the second induction hypothesis.

All that is left to do is to clarify the existence of the element $x'$ above. We will only point out which elements should be chosen as $x'$, the verification of the conditions is easily done. We can assume that $\dd(x,z) \geq 2$. 
\begin{itemize}
\item{$n=3$}
\begin{itemize}
\item $\dd(x,z)=2$, here we set $x'=z$.

\item $\dd(x,z)=3$, here we take $x' \inc z$, such that $u(xx',z) = 0$. The existence of such an $x'$ follows from applying Lemma~\ref{lemma:path} on a triangle containing $x,z$ and two elements incident with $x$ constructed by (U1).
\end{itemize}
\item{$n=4$}
\begin{itemize}
\item $\dd(x,z)=2$, here we set $x'=z$.
\item $\dd(x,z)=3$, let be $(x,a,b,z)$ the unique path of lenght 3 from $x$ to $z$. If $u(a,z)=0$, we let $x'$ be $b$. If this is not the case then let $c$ be an element incident with $x$ and such that $u(a,c)=0$. Next construct an element $d$ incident with $c$ such that $u(x,d)=0$. The last two constructions are possible by Lemma~\ref{lemma:pq0}. Finally $x'$ will be the projection of $d$ on $z$. Note that $x$ and $x'$ are equidistant due to Lemmas~\ref{lemma:oddquad} and~\ref{lemma:samedist}.
\item $\dd(x,z)=4$, if $x$ and $z$ are equidistant, we let $x'$ be $z$. Otherwise, using Lemma~\ref{lemma:choice}, we can construct a path $(x,a,b,c,z)$ such that $u(x,b) \geq u(b,z)$ and $u(a,c)=0$. Here we let $x'$ be the element $b$.

\end{itemize}

\item{$n=6$ and discrete}.
In this case the existence is guaranteed by the discreteness and Lemma~\ref{lemma:element}. 
\qed
\end{itemize}

\begin{cor}\label{cor:red}
If we reduce $v^{V(p,l) V(L,m) V(p',l') V(L',m')}$ to an expression of the form $v^{V(p'',l'') V(L'',m'')}$, then $l'' +m'' \leq l+m+l'+m'$.
\end{cor}
\proof
All the reductions in the proof of the above lemma use Corollary~\ref{cor:apart}, which does not increase the sum of the lengths of the translations.
\qed

\begin{lemma}\label{lemma:in3}
For each pair of valuations $v,w \in \Lambda(u)$ there is an ordinary $n$-gon $\Omega$ in $\Gamma$ which is non-folded for both $v$ and $w$.
\end{lemma}
\proof
Due to the previous lemma there exists a point $p$ and line $L \inc p$ in $\Gamma$, $l,m \in \R^+$ such that $w=v^{V(p,l) V(L,m)}$. Let $\Omega$ be an ordinary $n$-gon in $\Gamma$ containing $p$ and $L$ such that $\Omega$ is non-folded for $v$ (these exist because of Lemma~\ref{lemma:contain}). Because both $p$ and $L$ lie in $\Omega$, translations towards $p$ and $L$ produce valuations for which $\Omega$ remains non-folded. In particular this holds for $w=v^{V(p,l) V(L,m)}$.
\qed

\subsection{Building the affine apartment system}
We end by putting all the pieces together to form an affine apartment system. Let $\Lambda(u)$ be the set of points.  Remember that if $v \in \Lambda(u)$, then $\Lambda(u)=\Lambda(v)$. 

Let $\Omega$ be an ordinary $n$-gon of $\Gamma$. Consider the set $A(\Omega)$ of all the valuations in $\Lambda(u)$ for which this $n$-gon is non-folded. Suppose that two valuations $v_1$ and $v_2$ are in this set.  Lemma~\ref{lemma:2ok} tells us that there exists a flag $\{p,L\}$ in $\Gamma$ and $k,l \in \R^+$ such that $v_2=v_1^{V(p,k) V(L,l)}$. As $\Omega$ is non-folded for both $v_1$ and $v_2$, Lemma~\ref{lemma:convex} implies that there exists a flag $\{p',L'\}$ in $\Omega$ such that $v_2=v_1^{V(p',k) V(L',l)}$. We thus have that all the valuations in the set $A(\Omega)$ can be obtained out of each other by translating towards elements of $\Omega$. This is exactly the set of valuations which has been studied by Corollary~\ref{cor:apart}. In the reasoning behind this corollary it was seen that the valuations can be interpreted as points of $\A$. The sector with source $v \in \Lambda(u)$ and as direction the flag $\{p,L\}$ will be the set $\{v^{V(p,k)V(L,l)} |k,l \in \R^+ \}$.

This allows us to define a chart $f_{\Omega,v,p,L}$, for a $v\in \Lambda(u)$, and $\Omega$ a non-folded $n$-gon, containing a flag $\{p,L\}$ (the chart is defined such that a choosen fixed sector of $\A$ is mapped to the sector with source $v$ and direction $\{p,l\}$). Let $\cF$ be the collection of all these charts. Condition (A1) can now easily seen to be true. 



The second condition to check is (A2). Let $f=f_{\Omega,v,p,L}$ and $f'=f_{\Omega',v',p',L'}$ be two charts in $\cF$.  Let $X=f^{-1}(f'(\A))$. The points (or valuations) which are in the image of both charts, are those valuations for which both $\Omega$ and $\Omega'$ are non-folded. Let $v''$ be a valuation for which this is the case (if there is not such a $v''$, the condition (A2) is trivially satisfied). Lemma~\ref{lemma:convex} implies that $X$ is star convex for $f^{-1}(v'')$. Because $v''$ is arbitrary in $f(X)$, one obtains that $X$ is convex. That $X$ is also closed follows from the fact that translating changes the valuations continuously. 

Next thing we need to show is the existence of a $w \in W$ such that $f|_X = f' \circ w|_X$. Consider both $X$ and the similar set $X'=f'^{-1}(f(\A))$. In order to prove the existence of such a $w$ we need to prove that $X$ can be mapped on $X'$ by a $w \in W$. The map $\phi = f'^{-1} \circ f$ forms a bijective map from $X$ to $X'$. Let be $x_1$ and $x_2$ be elements of $X$. Then their images under $f$ are two valuations $v_1$ and $v_2$. Because they lie in the same apartment $A(\Omega)$, there is a flag $\{p,L\}$ in $\Omega$ and $k,l \in \R^+$ such that $v_2=v_1^{V(p,k) V(L,l)}$. But as these two valuations are also in $A(\Omega')$, we know by Lemma~\ref{lemma:convex} that there exists a flag $\{p',L'\}$ in $\Omega'$ such that $v_2=v_1^{V(p',k) V(L',l)}$. Because the lengths of the translations and the type of elements to which is translated do not change, it follows that $\phi$ is distance preserving and preserves the type of the directions at infinity of $\A$. This implies the existence of the needed $w$.

Condition (A3) is satisfied because of Lemma~\ref{lemma:in3}. 

Now, (A4) can be shown to be true as follows : suppose we have two sectors related to two flags $\{p,L\}$ and $\{q,M\}$ of $\Gamma$. These can be embedded in an ordinary $n$-gon $\Omega$, the apartment $A(\Omega)$ contains sectors with directions $\{p,L\}$ and $\{q,M\}$. This only leaves us to prove that two sectors related to the same flag always intersect in a subsector. This last assumption is true because if we have two ordinary $n$-gons $\Omega$ and $\Omega'$ containing $p$ and $L$, it follows from Corollary~\ref{cor:inc} that there exist $l,m \in \R^+$ such that for each $l' \geq l, m'\geq m$ the valuation $u^{V(p,l') V(L,m')}$ takes only the value zero in both $\Omega$ and $\Omega'$. The set of these valuations forms the desired subsector.

For (A5) we have three ordinary $n$-gons $\Omega$, $\Omega'$ and $\Omega''$, each pair sharing a path of length of $n$. From (U3) and (U4) we deduce that, if for a valuation $v \in \Lambda(u)$ the ordinary $n$-gon $\Omega$ is non-folded, then at least one of $\Omega'$ and $\Omega''$ is non-folded for $v$, too. This means that every point of $A(\Omega)$ belongs to $A(\Omega')$ or to $A(\Omega'')$, or to both. Since it is easy to see that the intersection of two apartments is closed, the sets $A(\Omega)\cap A(\Omega')$ and $A(\Omega)\cap A(\Omega'')$ are not disjoint, proving (A5).

We only still need to prove that the `distance' function $\dd$ defined on pairs of valuation by (A1), (A2) and (A3) is indeed a distance function. (For two valuations $v$ and $v^{V(p,k) V(L,l)}$, the distance between both will be the length of the third side of a triangle in an Euclidean plane where two sides have length $k$ and
$l$, and with the angle between both sides $\pi/n$.) However by rereading the proof in~\cite{parreau} of the equivalence of these definitions, one sees that showing the weaker inequality $\dd(u,v) \leq 2(\dd(u,w)+\dd(w,v))$ also suffices. This inequality is a direct consequence of Corollary~\ref{cor:red}.

So we conclude that the set of points $\Lambda(u)$, endowed with the set of apartments $\{A(\Omega)~|~ \Omega \mbox{ is an ordinary }n\mbox{-gon of } \Gamma \}$, forms a $2$-dimensional affine apartment system with the generalized $n$-gon $\Gamma$ at infinity. 

All that is left to show is that the construction of~\cite{kshvm} applied to the affine apartment system defined on $\Lambda(u)$ and the point defined by the valuation $u$, gives us back the valuation $u$ on $\Gamma$. One has to prove that, if $x$ and $y$ are adjacent, the corresponding sector-panels with source $u$ share a line segment of length $u(x,y)$. This follows from Corollary~\ref{cor:equiv} and the fact that, if $x$ and $y$ are adjacent, one has $\dd_r^{t,x}(y)=0$ if and only if $t \in [0,u(x,y)[$.

\section{Proof of the application (Theorem~\ref{theorem:app})}
Suppose we have given a projective plane $\Gamma$ and a real number $t\in \R^+ \backslash \{0\}$. Also suppose we are either given a valuation $u$, or two functions $\dd$ and $\angle$ satisfying the conditions listed in Theorem~\ref{theorem:app}. Use the identities $\dd(p,q)=t^{-u(p,q)}$ and $\angle{L,M} = \arcsin(t^{-u(L,M)})$ to reconstruct the other function(s).

It is easily seen that condition (U2) for valuations corresponds to condition (M2) and the part ``$\dd(p,q)=0 \Leftrightarrow p=q$'' of condition (M1).

If we have three points $p$, $q$ and $r$, then $$u(p,q)\geq \min(u(p,r), u(r,q)) \Leftrightarrow \dd(p,q)\leq \max(\dd(p,r), \dd(r,q)).$$ The left hand side is satisfied for a valuation because of (U3) and Lemma~\ref{lemma:path}; the right hand side is satisfied for a distance because of (M1). So condition (U3) for points on a line is equivalent with the inequality part of (M1).  

Condition (U1) for valuations is directly equivalent with conditions (M3) and (M4).

Also condition (U4) corresponds directly to the sine rule condition (M5). 

The only part that needs a closer look is how condition (U3) for valuations follows from conditions (M1) up to (M5) (and the already proven conditions (U1), (U2), (U3) for points on a line and (U4)). Let $L$, $M$ and $N$ be three lines through a point $p$. Using (U1), there exists two lines $Y$ and $Z$ through $p$ such that $u(Y,Z)=0$, and because we know (U1) and (U3) for points on a line hold, Lemma~\ref{lemma:pq0} also holds, so there exists a $q \inc Y$ and $r \inc Z$ with $u(p,q) =u(p,r) =0$. We now have for the line $qr$ that $\tau(p,qr)=0$ by (U4). (Note that $\tau$ is well-defined because (U4) holds.) 

Let $l$, $m$ and $n$ be the respective projections of $L$, $M$ and $N$ on the line $qr$ . Using (U4) we see that $u(L,M)=u(l,m)$, $u(M,N)=u(m,n)$ and $u(L,n)=u(l,n)$. So condition (U3) for the three lines $L$, $M$ and $N$ follows directly from the same condition (U3) for the three points $l$, $m$ and $n$.

\section{Some examples}
\subsection{$n=3$}
Here we rely on some results for the discrete case. The second author proved in \cite{Mal:93} that the notion of a projective plane with valuation is equivalent to one of a planar ternary ring with valuation. Moreover he also investigated in \cite{Mal:87} how the valuation behaves in planar ternary rings with extra algebraic properties (nearfields, quasifields, linear PTRs, etc.). In particular he proved the following result, th arguments of which can be copied verbatim in the non-discrete case.
\begin{prop}
A \emph{quasifield with valuation $v$}, which is an unary function with values in $\mathbb{Z} \cup \infty$ gives rise to a planar ternary ring with valuation (and thus also to a projective plane with valuation and an affine apartment system with a projective plane at infinity), if the following three conditions are fulfilled:
\begin{itemize}
\item[(V1)]
$v(a)= \infty$ if and only if $a=0$.
\item[(V2)]
If $v(a) < v(b)$, then $v(a+b) = v(a)$.
\item[(V3)]
$v(a_1 b - a_2 b) = v(a_1 - a_2) +v(b)$. 
\end{itemize}
\end{prop}   

We now construct such quasifields (again inspired by previous results of the second author in \cite{Mal:87}). Let $K_{+,.}$ be a field with a non-discrete valuation $v$ in the classical sense (which is in fact the above definition for quasifields applied to fields, so (V3) becomes $v(ab)=v(a) +v(b)$). 

\begin{remark}\rm
Notice that the classical affine apartment systems with a (Desarguesian) projective plane at infinity already appear here by taking those quasifields with valuation which are (skew)fields.
\end{remark}

Now let $\alpha$ be a field automorphism, with finite order, of $K$, preserving the valuation $v$. So $\alpha$ generates a finite group of automorphisms $G$. One can define the \emph{norm map} $n: K \rightarrow K : a \mapsto \prod_{\alpha' \in G} a$. Notice that $v(n(a))= |G| v(a)$. Let $\sigma$ be a map from the image of the norm map $n$ to $G$ such that $\sigma(1)$ is the unit element of $G$, and that $v(a)=v(b)$ implies $\sigma(n(a))=\sigma(n(b))$.

It follows that one can construct an Andr\'e quasifield $K_{+,\odot}$ by taking the elements of $K$ with the addition of the field and a new multiplication $\odot: K \times K \rightarrow K : (a,b) \mapsto a.b^{\sigma(n(a))}$. Moreover, we now show that this quasifield with the valuation $v$ forms a quasifield with valuation. We only have to verify (V3) for the new multiplication. First remark that $v(a\odot b) = v(a.b^{\sigma(n(a))})=v(a)+v(b^{\sigma(n(a))}) =v(a)+v(b)$. The last step holds because $\alpha$ and thus all elements of $G$ preserve $v$.

We now calculate $v(a_1 \odot b - a_2 \odot b)$. There are two possibilities that can occur.
\begin{itemize}
\item
$v(a_1)\neq v(a_2)$, suppose without loss of generality that $v(a_1) <v(a_2)$. Then
\begin{align}
v(a_1 \odot b - a_2 \odot b) &= v(a_1 \odot b) \\
&= v(a_1) + v(b) \\
&= v(a_1-a_2) +v(b), 
\end{align}
where the first step is true because $v(a_1 \odot b) = v(a_1) +v(b) < v(a_2)+v(b) = v(a_2)$, (V2), and $v(-1)=0$ (which easily follows from the definition of valuation).
\item
The other possibility is that $v(a_1)=v(a_2)$. Then
\begin{align}
v(a_1 \odot b - a_2 \odot b) &= v(a_1.b^{\sigma(n(a_1))} - a_2.b^{\sigma(n(a_2))})  \\
&= v((a_1-a_2).b^{\sigma(n(a_1))}) \\
&= v(a_1-a_2) +v(b^{\sigma(n(a_1))}) \\
&= v(a_1 -a_2) +v(b), 
\end{align}
where the second step holds because $v(a_1) = v(a_2)$ implies $\sigma(n(a_1)) = \sigma(n(a_2))$.
\end{itemize}
Combining both cases, we see that (V3) holds for the quasifield $K_{+,\odot}$ with valuation $v$. 

We end with some explicit examples of the above situation. Let $k$ be any field, let $M$ be a subset of $\mathbb{N}\backslash \{0\}$ generated multiplicatively by a certain set of primes. Now let $K$ be the field of rational functions in $t$, but allowing all rational powers $r/s$ of $t$ with $s\in M$. If $k(t)=f(t)/g(t)\in K$ with $f(t)$ and $g(t)$ polynomials (also allowing powers of the form above), we then set $v(k(t))$ to be minimal non-vanishing power of $t$ in $f(t)$ minus the minimal non-vanishing power of $t$ in $g(t)$. One verifies that $K$ together with $v$ forms a field with valuation. 

\begin{itemize}
\item Let $k$ be a finite field with characteristic $p$ and $M$ the set of integer powers of $p$. Then a suitable choice of $\alpha$ is the automorphism that maps $t^{\frac{r}{s}}$ to $(\frac{t^{1/s}}{1+t^{1/s}})^{r}$.
\item Now let $k$ be any field and $M$ generated by all the odd primes (so $M$ is the set of the odd non-negative integers). Now one can set $\alpha$ to be the automorphisms that maps $t^{\frac{r}{s}}$ to $(-t^{\frac{1}{s}})^r$.
\end{itemize}

All of these examples have a non-classical projective plane at infinity, but have classical residues. In addition they are locally finite when $k$ is finite. 

There are also examples where one can choose one residue completely freely. For a given planar ternary ring $R$, one can define a ``positively valuated ternary ring'' $R\{t\}$, similarly as in the discrete case, see \cite{Mal:88}. Indeed, one considers the power series $\sum_{n \in N}a_n t^n$ in $t$ where $N$ is a set of positive integer multiples of a certain rational number (for different power series, this number may be different) and $a_n\in R$ for $n\in N$. Since any finite number of such power series can be thought of as belonging to the same discrete version of this construction, the ternary operation can be copied from \cite{Mal:88}, and also the proof of the fact that we have a positively valuated ternary ring. Now, in completely the same way as in the discrete case, one constructs a projective plane with (non-discrete) valuation out of this.     The residue defined by this valuation is precisely the projective plane coordinatized by $R$. To the best of our knowledge, these are the first examples of such non-discrete apartment systems with an arbitrary residue.

\subsection{$n=4$}
The construction we will explain here is again inspired on an example for the discrete case by the second author in \cite{Mal:89}.  We will only sketch how the coordinatizing structure with valuation looks like. All proofs for the finite case still hold here (this is due to the fact than any finite number of elements in the coordinatizing structure can be `embedded' in a discrete case).  In particular, the reader can consult \cite{Mal:93} for explicit formulae to derive the valuation of the generalized quadrangle from the valuation of the coordinatizing structure.   

Consider the finite field $k=\mathsf{GF}(q)$ with $q=2^h$. Let $h_1$ and $h_2$ be two natural numbers such that $q-1$ and $-1+2^{1+h_1+h_2}$ are relatively prime (for example $h=3$, $h_1=1$ and $h_2=0$). For $i=1,2$, let $\theta_i$ be raising to the power $2^{h_i}$, which form automorphisms of this finite field. Now consider the field $K$ of Laurent series $\sum_{n \in N}a_n t^n$ in $t$ where $N$ is a set of integer multiples of a certain rational number, bounded below (again, for different Laurent series, this number may be different) and $a_n\in k$ for $n\in N$. There is a natural valuation on this field, defined by $v(\sum_{n \in N}a_n t^n)=m$ where $m$ is the smallest element of $N$ such that $a_m$ is non-zero (well defined by the boundedness below). One can extend $\theta_i$ for $i \in \{1,2\}$ to the field $K$ by 
\begin{equation}
(\sum_{n \in N}a_n t^n)^{\theta_i} = \sum_{n \in N}a_n^{\theta_i} t^n. 
\end{equation}
The coordinatizing structure is now given by:
\begin{align}
&Q_1(k,a,l,a') = (k^{\theta_1})^2 . a +a', \\
&Q_2(a,k,b,k') = a^{\theta_2}.k+k', 
\end{align}
wih $k,l,k',a,b,a' \in K$ and $v$ the natural valuation.

For more information about this example and coordinatizing structures see \cite{Mal:89}. One can show that this example defines a generalized quadrangle with valuation where both the quadrangle itself and its residue are non-classical. 

These are, to the best of our knowledge, the first explicitly defined examples of non-discrete $\R$-buildings of this nature.

\textbf{Address of the authors}\\
Department of Pure Mathematics and Computer Algebra,
Ghent University,\\
Krijgslaan 281-S22,\\
B-9000 Ghent,\\
BELGIUM\\
\texttt{\{kstruyve,hvm\}@cage.UGent.be}

\end{document}